\documentclass[preprint,review,12pt]{elsarticle}

\usepackage{amsmath,amssymb,amsthm}

\usepackage{caption}
\usepackage{multirow}

\usepackage{algorithm} 
\usepackage{algcompatible}
\usepackage{xcolor}
\usepackage{multicol}
\usepackage{tabularx,calc}

\newcommand{\EEA}{\end{eqnarray}}
\newcommand{\BEA}{\begin{eqnarray}}
\newcommand{\mycommenta}[1]{\small {\color{gray}\hspace{5mm}\%  #1   }}
\newcommand{\mycommentb}[1]{\small {\color{gray}\%  #1   }}
\newcommand{\comment}[1]{}

\usepackage{lineno,hyperref}
\modulolinenumbers[1]

\journal{XXX}

\begin{document}

\begin{frontmatter}

\title{Adaptive Error Covariances Estimation Methods for Ensemble Kalman Filters}

\author{Yicun Zhen\corref{cor}}
\cortext[cor]{Corresponding author.}
\ead{zhen@math.psu.edu}
\address{Department of Mathematics, The Pennsylvania State University, University Park, PA 16802}

\author{John Harlim}
\ead{jharlim@psu.edu}
\address{Department of Mathematics and Department of Meteorology, The Pennsylvania State University, University Park, PA 16802}

\begin{abstract}
This paper presents a computationally fast algorithm for estimating, both, the system and observation noise covariances of nonlinear dynamics, that can be used in an ensemble Kalman filtering framework. The new method is a modification of Belanger's recursive method, to avoid an expensive computational cost in inverting error covariance matrices of product of innovation processes of different lags when the number of observations becomes large. When we use only product of innovation processes up to one-lag, the computational cost is indeed comparable to a recently proposed method by Berry-Sauer's. However, our method is more flexible since it allows for using information from product of innovation processes of more than one-lag.

Extensive numerical comparisons between the proposed method and both the original Belanger's and Berry-Sauer's schemes are shown in various examples, ranging from low-dimensional linear and nonlinear systems of SDE's and 40-dimensional stochastically forced Lorenz-96 model. Our numerical results suggest that the proposed scheme is as accurate as the original Belanger's scheme on low-dimensional problems and has a wider range of more accurate estimates compared to Berry-Sauer's method on L-96 example. 

\end{abstract}


\end{frontmatter}


\section{Introduction}
Ensemble Kalman filters (EnKF) are ubiquitous in data assimilation of high-dimensional nonlinear problems \cite{Evensen1994,hunt:07}, particularly, it has received serious attentions in weather forecasting application \cite{houtekameretal:14}. While this practical method has been successful in applications, the accuracy of the resulting estimates is mostly determined by the prescription of both the system and observation error covariance statistics in the absence of model error. In fact, obtaining accurate covariance statistical estimates is still a challenging problem. Recent study \cite{ls:12} suggested that one should be cautious in interpreting the covariance estimates from approximate filtering methods, including EnKF.

Many numerical methods have been proposed to improve the filter covariance estimates. For example, a naive strategy is to directly inflate the covariance statistics with an empirically chosen inflation factor \cite{aa:99,ott:04,whwst:08}, while a second class of strategy is to use the current observations to adaptively inflate/tune the covariance statistics \cite{Mehra1970,Mehra1972,Belanger1974,DeeCohnDalcherGhil1985,anderson:07,lkm:09,BerrySauer2013,HarlimMahdiMajda2014}. 
The central contribution of this paper is on a new adaptive covariance estimation method that can be used with EnKF. In particular, the new method is motivated by two recently developed EnKF-based covariance estimation methods \cite{BerrySauer2013,HarlimMahdiMajda2014} that generalized, respectively, two competing methods that were introduced in early 70's for covariance estimation of stationary \cite{Mehra1970,Mehra1972} and time-varying linear problems \cite{Belanger1974}. All of these methods share a fundamental similarity, that is, they use the information from lag correlation of innovation processes (or prior forecast errors in the observation space) to approximate $Q$ and $R$. While Berry-Sauer's estimation method \cite{BerrySauer2013} is computationally fast, their method used only product of innovation processes up to one-lag, $L=1$, which by design, restricts its applicability especially when observations are sparse. Furthermore, the accuracy of the resulting estimates can be sensitive to the choice of a nuisance parameter, as we will show later. On the other hand, the covariance estimation method of \cite{HarlimMahdiMajda2014}, which is essentially Belanger's method implemented with EnKF \cite{Belanger1974}, allows for using product of innovation processes of more than one-lag, $L>1$,  and produces estimates that converge faster. However, this method is computationally expensive when the number of observations, $m$, become large because it involves inversion of $m^2\times m^2$ matrix for $L+1$ times in each time step. This computational cost can be reduced from $\mathcal{O}((N_pm^6+N_p^2m^2+N_p^3)L)$ to $\mathcal{O}((N_pm^3+N_p^2m^2+N_p^3)L)$ by observing a special tensor structure of the covariance matrix, where $N_p$ is the number of parameters for $Q$ and $R$ \cite{DeeCohnDalcherGhil1985}. However, this reduced computational cost may still be unaffordable when the number of observations is large. In this paper, we present a new estimation method with computational cost $\mathcal{O}(m^2N_p^2(L+1))$, that is faster than the Belanger's scheme when $N_p\ll m$, and more robust than Berry-Sauer's scheme \cite{BerrySauer2013} when observations are sparse, and produces accurate estimates. We will provide the computational costs when these methods are used in tandem with ETKF (\cite{hunt:07} version). We will refer to this new method as the modified Belanger's method since it is motivated by the original Belanger's formulation. We shall see that some aspect of the proposed method is also motivated by Berry-Sauer's scheme. 

The remainder of this paper is organized as follows: In Section~2, we will discuss the key idea behind the scheme and point out the similarities and differences between the proposed method with the existing schemes. This discussion is supported by two Appendices: Appendix A describes the detailed derivation and Appendix B provides pseudo-codes for an easy access to implementation. In particular, we will implement this scheme with ETKF that was used in \cite{HarlimMahdiMajda2014}, which was formulated by \cite{hunt:07}. In Section~3 we show numerical results, comparing the proposed method with the original Belanger's scheme implemented via ETKF \cite{HarlimMahdiMajda2014} and Berry-Sauer's scheme \cite{BerrySauer2013, Belanger1974} on various test models ranging from low-dimensional linear, low-dimensional nonlinear, and moderately higher-dimensional nonlinear stochastic filtering problems. We also include an experiment with LETKF \cite{hunt:07} to demonstrate the potential for high-dimensional applications. We close the paper with a short summary and discussion in section~4.

\section{Methodology}
In this section, we present the mathematical formulation of a new method for estimating noise covariance matrices. Since the new method is motivated by earlier methods \cite{Belanger1974,BerrySauer2013,OdelsonRajamaniRawlings2006}, we provide a short review of these existing methods and point out the similarities and differences. To simplify the presentation in this section, we assume that the observations are taken at every integration time step. We also provide the formulation for general case in which the observations are taken at every $N\geq1$ integration time steps in Appendix A and the pseudo-algorithms in Appendix B.  

\subsection{Mathematical Formulation}
We consider the following discrete-time linear filtering problem:
\BEA
x_{j}&=&F_{j-1}x_{j-1}+\Gamma w_{j-1},\label{model}\\
y^o_{j}&=&H_{j}x_j+\xi_j,\label{obsmodel}
\EEA
where $x_j\in\mathbb{R}^n$ denotes the hidden state variable at time $t_j$. In \eqref{model}, we assume that the dynamical model is driven by $\ell-$dimensional independent identically distributed (i.i.d.)~Gaussian noises $w_{j-1}$ with known $\Gamma\in\mathbb{R}^{n\times\ell}$ and an unknown time-independent system error covariance matrix, $Q\in\mathbb{R}^{q\times q}$. The observations $y^o_{j}\in\mathbb{R}^m$ in \eqref{obsmodel} are corrupted by $m-$dimensional Gaussian i.i.d~noises, $\xi_j$, with an unknown time-independent observation error covariance matrix, $R\in\mathbb{R}^{m\times m}$. Throughout this paper, we assume that $F_j$, $H_j$ are linear operators independent of the system and observation noises. Our goal is to estimate the hidden state variable $x_j$ and the noise covariance matrices, $Q$ and $R$, given only noisy observations, $y^o_j$.

If the noise covariances $Q$ and $R$ are known, an unbiased optimal estimation (in the sense of minimum variance) for state $x_j$ can be obtained by Kalman filter formula \cite{KalmanBucy1961}. 
In particular, the Kalman filter formula recursively updates the prior mean and covariance estimates, $x_j^f$ and $B_j^f$, respectively, of a conditional distribution $p(x_j |y^o_i, i\leq j-1)$, to the posterior mean and covariance estimates, $x_j^a$ and $B_j^a$ of a new conditional distribution $p(x_j |y^o_i, i\leq j)$, incorporating the observation $y^o_j$ at time $t_j$ (see Appendix A for the detailed of the Kalman filter formula, cf.~\eqref{KF_N}). 

Following the formulation in \cite{Belanger1974}, we will use the {\it innovation process}, which is defined as the error of the prior mean estimates in the observation subspace,
\BEA
v_{j}:=y^o_j-H_jx^f_j = H_j \Delta x^f_j + \xi_j ,\label{innovation}
\EEA
to estimate $Q$ and $R$. In particular, Belanger's formulation relies on the description of prior forecast error with the following recursive equation \cite{Belanger1974},
\BEA
\Delta x^f_j&:=&x_j-x^f_j = \prod_{i=0}^{j-1}\mathcal{U}_i\Delta x^f_0-\mathcal{G}^{\xi}_{(j)}+\mathcal{F}^{w}_{(j)},\label{Deltaxf}
\EEA
where the explicit definition of the terms $\mathcal{U}_{j}$,  $\mathcal{G}^{\xi}_{(j)}$, and $\mathcal{F}^{w}_{(j)}$, 
are presented in appendix A (cf. \eqref{eq28}, \eqref{gxij}, \eqref{fwj}). We remark that $\mathcal{G}^{\xi}_{(j)}$ and $\mathcal{F}^w_{(j)}$ are linear combinations of, respectively, $\xi$ and $w$ of previous time steps with coefficients involving  $F_{j-i}$,  $H_{j-i}$ and the Kalman gain matrix $K_{j-i}$. Hence the coefficients of $\xi$ in $\mathcal{G}^{\xi}_{(j)}$ and $w$ in $\mathcal{F}^w_{(j)}$ are independent of the realization of  $\xi$ and $w$ as long as the choice of Kalman gain matrix $K_{j-i}$ at each time step is independent of $\xi$ and $w$.  However, this requirement can not be satisfied in any adaptive filter method which estimates $Q$ and $R$ sequentially since this estimation depends on the realization of $\xi$ and $w$. In this situation, $\mathcal{G}^{\xi}_{(j)}$ and $\mathcal{F}^{w}_{(j)}$ are no longer linear functions of $\xi$ and $w$. But they still can be viewed as linear combinations of $\xi$ and $w$ whose coefficients also involve $\xi$ and $w$.

Taking the expectation of the product of \eqref{innovation} at different lags with respect to realizations of $\xi$ and $w$, we have:
\BEA
\mathbb{E}[v_jv_{j-l}^{\top}]&=&H_j\mathbb{E}[\Delta x^f_j(\Delta x^f_{j-l})^{\top}]H_{j-l}^{\top}+H_j\mathbb{E}[\Delta x^f_j\xi_{j-l}^{\top}]\nonumber\\
&=&H_j{\Big \{}\mathbb{E}[\mathcal{G}^{\xi}_{(j)}(\mathcal{G}^{\xi}_{(j-l)})^{\top}]-\mathbb{E}[\mathcal{G}^{\xi}_{(j)}(\mathcal{F}^w_{(j-l)})^{\top}]\nonumber\\
&&-\mathbb{E}[\mathcal{F}^w_{(j)}(\mathcal{G}^{\xi}_{(j-l)})^{\top}]+\mathbb{E}[\mathcal{F}^w_{(j)}(\mathcal{F}^w_{(j-l)})^{\top}]{\Big \}}H_{j-l}^{\top}\nonumber\\
&&-H_j\mathbb{E}[\mathcal{G}^{\xi}_{(j)}\xi_{j-l}^{\top}]+H_{j}\mathbb{E}[\mathcal{F}^{w}_{(j)}\xi_{j}^{\top}] + \text{remainder terms},\label{eq5}
\EEA
where all the terms that contain $\prod_{i=0}^{j-1}\mathcal{U}_i$ are absorbed in the remainder terms. When the filter is uniformly asymptotically stable, these remainder terms decay to zero at an exponential rate \cite{Belanger1974,OdelsonRajamaniRawlings2006}. If the Kalman gain matrices $K_{j-i}$ are computed independently of $\xi$ and $w$, the cross terms in \eqref{eq5} can be eliminated and the right hand side of \eqref{eq5} is a linear function of $Q$ and $R$. Let's parameterize $Q$ and $R$ with $N_Q$ and $N_R$ parameters using some prescribed basis $Q_s$ and $R_s$, respectively, such that:
\BEA
Q=\displaystyle\sum_{s=1}^{N_Q}\alpha_sQ_s, \quad R=\displaystyle\sum_{s=1}^{N_R}\beta_sR_s.\nonumber
\EEA
Then, we can approximate \eqref{eq5} as follows,
\BEA
\mathbb{E}[v_jv_{j-l}^{\top}]\approx\sum_{s=1}^{N_Q}\alpha_s\mathcal{H}^{(Q)}_{j,l,s}+\sum_{s=1}^{N_R}\beta_s\mathcal{H}^{(R)}_{j,l,s},\label{eq6}
\EEA
where the explicit definition of $\mathcal{H}^{(Q)}_{j,l,s}$ and $\mathcal{H}^{(R)}_{j,l,s}$ are provided in Appendix A (cf.~\eqref{HQjls}, \eqref{HRjls}) and the approximation here is due to neglecting the remainder terms in \eqref{eq5}. 

In general situation when the Kalman gain matrix depends on the realization of $\xi$ and $w$, the approximation in \eqref{eq6} is more than just due to neglecting the remainder terms. This is because the cross terms in \eqref{eq5} are not zero and the right hand terms in \eqref{eq5} are no longer linear functions about $Q$ and $R$. Now, let $\mathcal{E}_{j,l}$ be defined as the error of the approximation in \eqref{eq6}. Taking average of these errors over time,
\BEA
\frac{1}{J-L}\sum_{j=L+1}^{J}\mathcal{E}_{j,l} &:=& \frac{1}{J-L} \sum_{j=L+1}^J \Big(\mathbb{E}[v_{j}v_{j-l}^{\top}] - (\sum_{s=1}^{N_Q}\alpha_s\mathcal{H}^{(Q)}_{j,l,s}+\sum_{s=1}^{N_R}\beta_s\mathcal{H}^{(R)}_{j,l,s})\Big)\nonumber \\
&=& \frac{1}{J-L} \sum_{j=L+1}^J v_jv_{j-l}^{\top} - \frac{1}{J-L} \sum_{j=L+1}^J \delta_{j,l}  \nonumber\\ & & -\frac{1}{J-L} \sum_{j=L+1}^J  (\sum_{s=1}^{N_Q}\alpha_s\mathcal{H}^{(Q)}_{j,l,s}+\sum_{s=1}^{N_R}\beta_s\mathcal{H}^{(R)}_{j,l,s})\label{eq7}
\EEA
where $\delta_{j,l}:=v_{j}v_{j-l}^{\top}-\mathbb{E}[v_{j}v_{j-l}^{\top}]$. Suppose that the estimates for $Q_j$ and $R_j$ that are used to update the Kalman gain $K_{j}$ equilibrate to some constant values. Then $v_{j}$ is a stationary, Gaussian process with an exponential decaying correlation function $\|\mathbb{E}(v_iv_j^{\top})\|_{max}=\mathcal{O}(\exp(-\alpha |i-j|))$, where the order of magnitude is with respect to $|i-j|$. Then, it is not difficult to show that, 
\BEA
\mathbb{E}\Big[(\frac{1}{J-L}\displaystyle\sum_{j=L+1}^{J}\delta_{j,l})(\frac{1}{J-L}\displaystyle\sum_{i=L+1}^{J}\delta_{i,l})^\top\Big] = \mathcal{O}\Big(\frac{1}{J-L}\Big). \nonumber
\EEA
Secondly, the terms $\mathcal{E}_{j,l}$ are automatically eliminated when the Kalman gain matrices are updated with constant $Q_j$ and $R_j$. Hence it is reasonable to believe that the variance and bias of $\frac{1}{J-L}\displaystyle\sum_{j=L+1}^{J}\mathcal{E}_{j,l}$ are small when $Q_j$ and $R_j$ vary slowly as functions of time, resulting to the following approximation,
\BEA
\frac{1}{J-L} \sum_{j=L+1}^J v_jv_{j-l}^{\top} \approx \frac{1}{J-L}  \sum_{j=L+1}^J  (\sum_{s=1}^{N_Q}\alpha_s\mathcal{H}^{(Q)}_{j,l,s}+\sum_{s=1}^{N_R}\beta_s\mathcal{H}^{(R)}_{j,l,s}).\label{averageeqn}
\EEA

\subsubsection{Modified Belanger's Method}
Based on this observation, we propose the following method to estimate $\alpha_s$ and $\beta_s$:
\begin{itemize}
\item At the time point $t_J$, solve \eqref{averageeqn} for $\alpha_{J}:=(\alpha_{1,J},\ldots,\alpha_{N_Q,J})$ and $\beta_{J} := (\beta_{1,J},\ldots,\beta_{N_R,J})$ using a least square method,
\BEA
\min_{\alpha_{J},\beta_{J}} \sum_{l=0}^L\Big\|\sum_{j=L+1}^Jv_jv_{j-l}^{\top}-{\Big \{}\sum_{s=1}^{N_Q}\alpha_{s,J}(\sum_{j=L+1}^{J}\mathcal{H}^{(Q)}_{j,l,s})+\sum_{s=1}^{N_R}\beta_{s,J}(\sum_{j=L+1}^{J}\mathcal{H}^{(R)}_{j,l,s}){\Big \}}\Big\|\label{newmethod}
\EEA
where the new subscript $J$ denotes the time $t_{J}$ and the norm is Frobenius norm.
\item Relax the estimates of $\alpha$ and $\beta$ back to the estimates of the previous time step by the following running average:
\BEA
\alpha_{J}=\alpha_{J-1}+\frac{1}{\tau}(\alpha_{J}-\alpha_{J-1}),\label{Relx1}\\
\beta_{J}=\beta_{J-1}+\frac{1}{\tau}(\beta_{J}-\beta_{J-1}),\label{Relx2}
\EEA
where $\tau\geq 1$ is a nuisance parameter that will be empirically chosen.

\end{itemize}
We provide pseudo-algorithm \ref{CovEst_k} to guide the detail implementation of this method in Appendix B.


The idea of the first step of this method is to match the observed data $v_{j}v_{j-l}^{\top}$ with the estimation of its statistical mean in a least square manner since we have no access to $\mathbb{E}[v_{j}v_{j-l}^{\top}]$. A closely related method has been proposed by
\cite{OdelsonRajamaniRawlings2006}. Specifically, their method can be thought of as the stationary version of this method. They assume stationarity of $v_jv_{j-l}^{\top}$ and use the ergodicity property to approximate 
$\mathbb{E}[v_{j}v_{j-l}^{\top}] \approx\frac{1}{J-L}\displaystyle\sum_{j=L+1}^{J}v_jv_{j-l}^{\top}$ in \eqref{eq6}, and solve a different least square function, 
\BEA
\min_{\alpha_{J},\beta_{J}} \sum_{l=0}^{L}\Big\|\frac{1}{J-L}\sum_{j=L+1}^{J}v_jv_{j-l}^{\top}-\sum_{s=1}^{N_Q}\alpha_{s,J}\mathcal{H}^{(Q)}_{l,s}-\sum_{s=1}^{N_R}\beta_{s,J}\mathcal{H}^{(R)}_{l,s}\Big\|,\label{eq10}
\EEA
where the norm is also taken to be the Frobenius norm. An important difference between this method and the modified Belanger's method in \eqref{newmethod} is that the operators $\mathcal{H}^{(Q)}, \mathcal{H}^{(R)}$ in \eqref{eq10} are time independent and estimated separately assuming that the underlying processes are stationary \cite{OdelsonRajamaniRawlings2006}. However for any time-varying system their method is not applicable since we can't determine $\mathcal{H}^{(Q)}$ and $\mathcal{H}^{(R)}$. From this point of view, the modified Belanger's method is advantageous since it handles non-stationary cases by design.

The relaxation step (or the running average \eqref{Relx1} and \eqref{Relx2}), which was proposed in \cite{BerrySauer2013} for different reason, is included in this method to reduce the dependence of the Kalman gain matrices on $\xi$ and $w$ which subsequently reduce bias terms that can possibly occur when $\mathcal{E}_{j,l}$ and $\delta_{j,l}$ are large in \eqref{eq7}. From the practical point of view, this relaxation step also increases the stability (or reduce the covariance) of the estimate. 

\subsubsection{The original Belanger's Method}
 
The original Belanger's approach \cite{Belanger1974} mitigates the non-stationarity with a secondary Kalman filter, treating eqn~\eqref{eq6} as an observation operator for $Q$ and $R$. In particular, he assumes that the error in the approximation in \eqref{eq6} is white and normally distributed with mean zero and covariance $W_{j,l}\in\mathbb{R}^{m^2\times m^2}$, obtained through a Gaussian approximation. Subsequently, a secondary Kalman filter is employed to solve the following optimization problem, 
\BEA
\min_{\alpha_{J},\beta_{J}}  {\bf Q}(\alpha_J,\beta_J)+\sum_{j=L+1, l=0}^{J,L}\Big\|{\bf vec}(v_{j}v_{j-l}^{\top}-\sum_{s=1}^{N_Q}\alpha_{s,J}\mathcal{H}^{(Q)}_{j,l,s}-\sum_{s=1}^{N_R}\beta_{s,J}\mathcal{H}^{(R)}_{j,l,s})\Big\|^2_{W_{j,l}^{-1}},\label{belanger}
\EEA
where function ${\bf Q}$ is any arbitrarily chosen positive definite quadratic function about $\alpha_J =(\alpha_{1,J},...,\alpha_{N_Q,J})$ and $\beta_{J} = (\beta_{1,J},...,\beta_{N_R,J})$ and  the operator ${\bf vec}$ denotes the vectorization of a matrix.
Hence Belanger's method requires the inversion of $W_{j,l}$. 

As mentioned in the introduction, this step of matrix inversion can be computed in $\mathcal{O}((N_pm^3+N_p^2m^2+N_p^3)L)$ operations, realizing a special structure in $W$ \cite{DeeCohnDalcherGhil1985}, 
where we define $N_p=N_Q+N_R$ as the total number of parameters. On the other hand, the cost of minimizing \eqref{newmethod} in the modified Belanger's method requires only $\mathcal{O}(N_p^2m^2(L+1))$. In the context of ETKF, both the modified and original Belanger's methods construct matrices (such as $\mathcal{H}^{(Q)}$ and $\mathcal{H}^{(R)}$) to be used in the secondary filter and the computational costs are $\mathcal{O}(nmN_eN_p(L+1))$. Assuming that the observation error covariance $R$ is diagonal, the cost of  the primary filter (ETKF) is $\mathcal{O}(mN_e^2+N_e^3+nN_e^2)$ which is small compared to the total cost of the secondary filters; $\mathcal{O}((nmN_eN_p+N_p m^3+N_p^2m^2+N_p^3)L)$ for the original Belanger's method and $\mathcal{O}((nmN_eN_p+N_p^2m^2)(L+1))$ for the modified Belanger's method. Setting the number of parameters to be less than the observations, $N_p\ll m$, then the second term $\mathcal{O}(N_p m^3L)$ dominates the computational cost of the original Belanger's scheme if $m^2 \gg nN_e$. On the other hand, the second term in the computational estimate for the modified Belanger's scheme, $\mathcal{O}(N_p^2m^2L)$, is always less than the first term, $\mathcal{O}(nmN_eN_pL)$ whenever $N_pm\ll nN_e$. 

Therefore, to achieve the most efficient computational cost with the modified Belanger's method, one should parameterize $Q$ and $R$ with a total number of parameters $N_p=N_Q+N_R$ that is less than the total number of observations, $m$, and choose an ensemble size, $N_e$, appropriately to satisfy  $N_pm\ll nN_e$.

\subsubsection{Berry-Sauer's Method} 

Finally, we will compare our method to a computationally faster, related, method that was recently proposed by \cite{BerrySauer2013}. This method is a variant of Mehra's method \cite{Mehra1970} for non-stationary (nonlinear) processes. In particular, their method is based on a formulation of the expectation of the product of innovations, $\mathbb{E}[v_jv_{j-l}^\top]$, that does not use the recursive identity in \eqref{Deltaxf}; resulting to a much simpler set of equations for zero-lag and one-lag products,
\BEA
\mathbb{E}[v_jv_{j}^\top] &=& H_j B_j^f H_j^\top + R, \label{RQ} \\
\mathbb{E}[v_{j+1}v_{j}^\top] + H_{j+1}F_jK_j\mathbb{E}[v_jv_{j}^\top] &=& H_{j+1} F_j (F_{j-1}B^a_{j-1}F_{j-1}^\top + \Gamma Q\Gamma^\top)H_j^\top, \nonumber
\EEA  
where $B_j^f := \mathbb{E}[(x_j-x_j^f)(x_j-x_j^f)^\top]$ and $B_j^a := \mathbb{E}[(x_j-x_j^a)(x_j-x_j^a)^\top]$. The goal is to obtain $Q$ and $R$, which are statistics of stationary processes $w$ and $\xi$, respectively, from eqns~\eqref{RQ} that involve statistics of non-stationary processes, namely, $B_j^f, B_j^a, K_j$. Furthermore, we have no access to the statistical quantities $\mathbb{E}[v_jv_{j}^\top]$ and $\mathbb{E}[v_{j+1}v_{j}^\top]$.

Their method consists of two steps: First, they estimate $R$ and $Q$ separately by solving,
\BEA
 R = v_jv_{j}^{\top} - H_j \tilde{B}_j^f H_j^\top \label{solveR},
 \EEA
 and, 
 \BEA
\min_{\alpha_j}  \Big\| v_{j+1}v_{j}^\top + H_{j+1}F_jK_jv_jv_{j}^\top - \sum_{s=1}^{N_Q} H_{j+1} F_j (F_{j-1}\tilde{B}^a_{j-1}F_{j-1}^\top + \alpha_{s,j} \Gamma Q_s \Gamma^\top)H_j^\top  \Big\|\label{solveQ},
\EEA
where they use parameterization $Q=\displaystyle\sum_{s=1}^{N_Q}\alpha_{s,j}Q_s$ with the Frobernius norm in the cost function above, and define $\alpha_j :=(\alpha_{1,j}, \ldots, \alpha_{N_Q,j})$.
Additionally, they use the covariance estimates from the Kalman filter formula as an approximation to the theoretical covariance estimates, $\tilde{B}_j^f \approx B_j^f$ and $\tilde{B}_j^a \approx B_j^a$. The second step in their method is to apply the moving average for the estimates of $R$ and $Q$ as in \eqref{Relx1} to reduce the sensitivity (or variance) of the sequential estimates in \eqref{solveR} and \eqref{solveQ} which depend on realization of noises $w$ and $\xi$ through $v_jv_j^\top$ and $v_{j+1}v_{j}^\top$. In some sense, this method is analogous to the original Belanger's method except that it uses moving average as the secondary filter rather than Kalman filter. This cheaper secondary filter is also adopted in the modified Belanger's scheme, replacing the expensive secondary Kalman filter.
Thus, the computational cost of the modified Belanger's scheme is comparable to that of Berry-Sauer's method when the innovation lag, $L=1$.

While this method is computationally the cheapest one relative to the other methods discussed in this paper, it is not clear that the approximations in \eqref{solveR} and \eqref{solveQ} would induce bias in a long run. On the other hand, it is clear that bias indeed exists at any specific time point when the filter covariance estimates, $\tilde{B}_j^f \approx B_j^f$ and $\tilde{B}_j^a \approx B_j^a$, are sub-optimal since 
the endogenous variables in \eqref{solveR} and \eqref{solveQ} depend on these covariances. On the other hand, the endogenous variables of both the modified and original Belanger's methods in \eqref{newmethod} and \eqref{belanger}, respectively, involve only $v_{j+l}v_{j}^\top$. We will numerically show in Section~3 that this method is more sensitive to the choice of the nuisance parameter $\tau$ in the relaxation step compared to the proposed modified Belanger's method and it requires large enough $\tau$ to reduce the variance of the error of the estimates which subsequently slows down the convergence of the estimates. One practical limitation of this method is that it only uses information up to one-lag of innovation statistics which restricts its application when observations are sparse, as we will demonstrate in Section~3 below.


\subsection{Covariance estimation with ensemble Kalman filters}
The formulation of the  covariance estimation method described in the previous section depends on knowing the linear operators $F_j$ and $H_j$. For nonlinear problems, however, we typically know the nonlinear forward operator $f_{j}$ and observational operator $h_{j}$ instead of the linear operators. Here we will approximate the corresponding linear operators with an ensemble of solutions as in \cite{BerrySauer2013,HarlimMahdiMajda2014} which is natural in the ensemble Kalman filtering setting. 

More specifically, let $\{X^{a,i}_{j}\}_{i=1}^{N_e}$ be an ensemble of the posterior estimates of $x_{j}$ and $X^{df,i}_{j+1}=f_{j}(X^{a,i}_{j})$ be the $i-$ensemble member of deterministic forecasts from initial condition $X^{a,i}_{j}$. Let $U_{j}^{a}$ and $U_{j+1}^{df}$ be matrices whose $i$th columns are the ensemble perturbations of $X^{a,i}_{j}$ and $X^{df,i}_{j+1}$ from their ensemble mean, respectively, then the linear forward operator $F_{j}$ can be approximated by
\BEA
F_{j}\approx U^{df}_{j}(U^{a}_{j})^{\dagger}.
\EEA
where $\dagger$ denotes the matrix pseudo-inverse. 

Similarly, let the prior ensemble be defined as $\{X^{f,i}_{j}\}_{i=1}^{N_e}$ of the stochastic system (these ensemble are Gaussian random samples with mean $x_j^f$ and covariance $P_j^f$, cf.~\eqref{pjf}). Also, let $Y^{f,i}_{j}=h_{j}(X^{f,i}_{j})$ be the $i$th ensemble member that is (nonlinearly) projected to the observational space. Define $U^{f}_{j}$ and $V^{f}_{j}$ as matrices whose $i$th columns consist of the ensemble perturbation of $X^{f,i}_{j}$ and $Y^{f,i}_{j}$ from its ensemble mean, respectively. Then the linear observation operator can be approximated by:
\BEA
H_{j}\approx V^{f}_{j}(U^{f}_{j})^{\dagger}.
\EEA
Therefore, we can directly apply the new covariance estimation method discussed in Section~2.1.1 to update $Q$ and $R$ after each ensemble Kalman filter analysis step that updates the estimates of the state variables, $x_j$. 

A second adjustment for using EnKF is in the generation of a new ensemble forecast that accounts for the stochastic terms once $Q$ is updated. In particular, since 
\BEA
P^f_{j+1} = \frac{1}{N_e-1}U_{j+1}^{df}(U_{j+1}^{df})^\top + \Gamma Q\Gamma^\top,\label{pjf}
\EEA  
where $Q$ is obtained from the covariance estimation method, then we need to construct matrix $U_{j+1}^{f}$ such that $P^f_{j+1}  = (N_e-1)^{-1}U_{j+1}^{f} (U_{j+1}^{f})^\top$. In the numerical experiments below, we use the method proposed in \cite{HarlimMahdiMajda2014} which 
defines each column of $U_{j+1}^{f}$ as an $n-$dimensional normally distributed random vector with mean zero and covariance $P^f_{j+1}$. In Algorithm~\ref{ETKF_k} in the Appendix B, we will show the details of these adjustments for an implementation with Ensemble Transform Kalman Filter \cite{Bishop.et.al2001,hunt:07}. 

For the implementation with LETKF in Section~3.4, we approximate $F_j$ locally and we incorporate $\Gamma Q\Gamma^\top$ globally in the primary filter, even if the estimates of $Q$ are obtained independently from the local secondary filter. In our implementation, we define the global estimate of $Q$ to be the spatial average of local estimates of $Q$ and we construct the global $U^f_{j+1}$ by adding Gaussian random noises with mean zero and covariance $\Gamma Q\Gamma^\top$ to each column of the global prior ensemble perturbation, $U^{df}_{j+1}$. Also, in each local primary filter analysis (LETKF) we use the global estimate of $R$ (restricted to each local region), which we define as the spatial average of the local estimates of $R$ that are obtained by independent local secondary filter. 

\section{Numerical Simulations}
In this section, we will numerically compare the three methods discussed above, in particular, the proposed modified Belanger's scheme, the original Belanger's scheme, and Berry-Sauer's scheme on three simple models: 1) a two-dimensional linear model, 2) a low-dimensional nonlinear triad model, and 3) a stochastically forced 40-dimensional Lorenz-96 model. Finally, we will also compare results with LETKF, applying it to the Lorenz-96 model.

\subsection{A linear example}
In this section, we consider a two-dimensional linear filtering problem in \eqref{model} and \eqref{obsmodel} with parameters as in \cite{Belanger1974,Mehra1970}, 
\BEA
F_{j}=\begin{pmatrix}
0.75 & -1.74\\
0.09 & 0.91
\end{pmatrix}, \Gamma=\begin{pmatrix}
1 & 0.4\\
0.1 & 1
\end{pmatrix}, Q=\mathcal{I}_2, R= 0.5\mathcal{I}_m.\nonumber
\EEA
In the simulations below, we will consider full and sparse observations. Full observations correspond to $H_j=\mathcal{I}_2$ and a $2\times 2$ diagonal $R = 0.5\mathcal{I}_2$. On the other hand, sparse observations correspond to $H_j=[1,0]$ and a scalar $R=0.5$. Since $Q$ and $R$ are both diagonal matrices we naturally parameterize them by their diagonal components, 
\BEA
Q = \sum_{s=1}^{N_Q} \alpha_s  Q_s, \quad R = \sum_{s=1}^{N_R} \beta_s  R_s, \label{param}
\EEA
where in the case of full observations, $N_Q=N_R=2$ and $Q_s = R_s \in\mathbb{R}^{2\times 2}$ are matrices with one on the $s$th diagonal component and zero everywhere else. In the case of sparse observations $R_s = 1$ but $Q_s$ remains the same. 
In the numerical simulations below, we will show results for 10000 and 50000 assimilation cycles for the cases of  full and sparse observation, respectively. 

\begin{figure}
\includegraphics[width=0.5\textwidth]{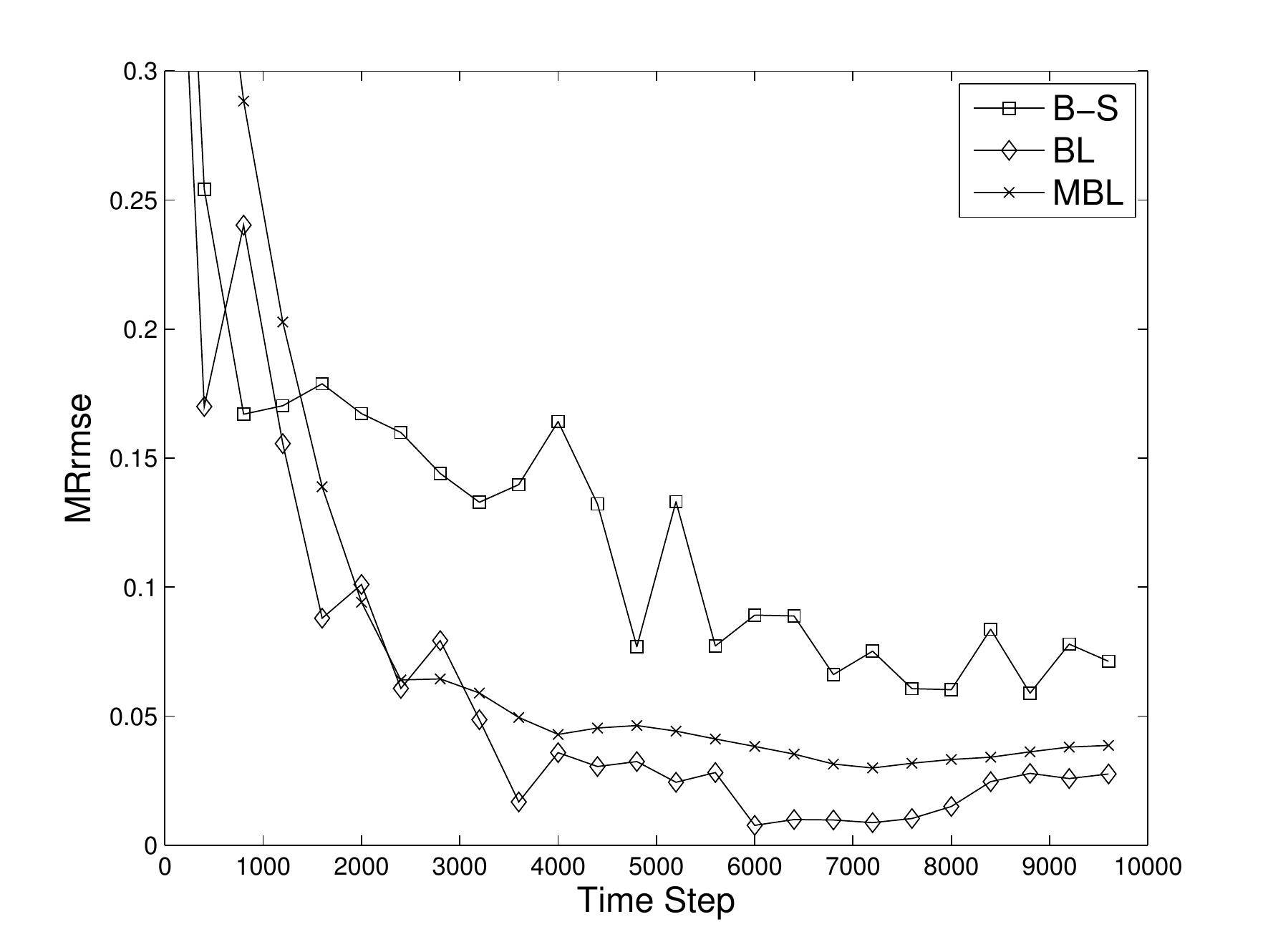}
\includegraphics[width=0.5\textwidth]{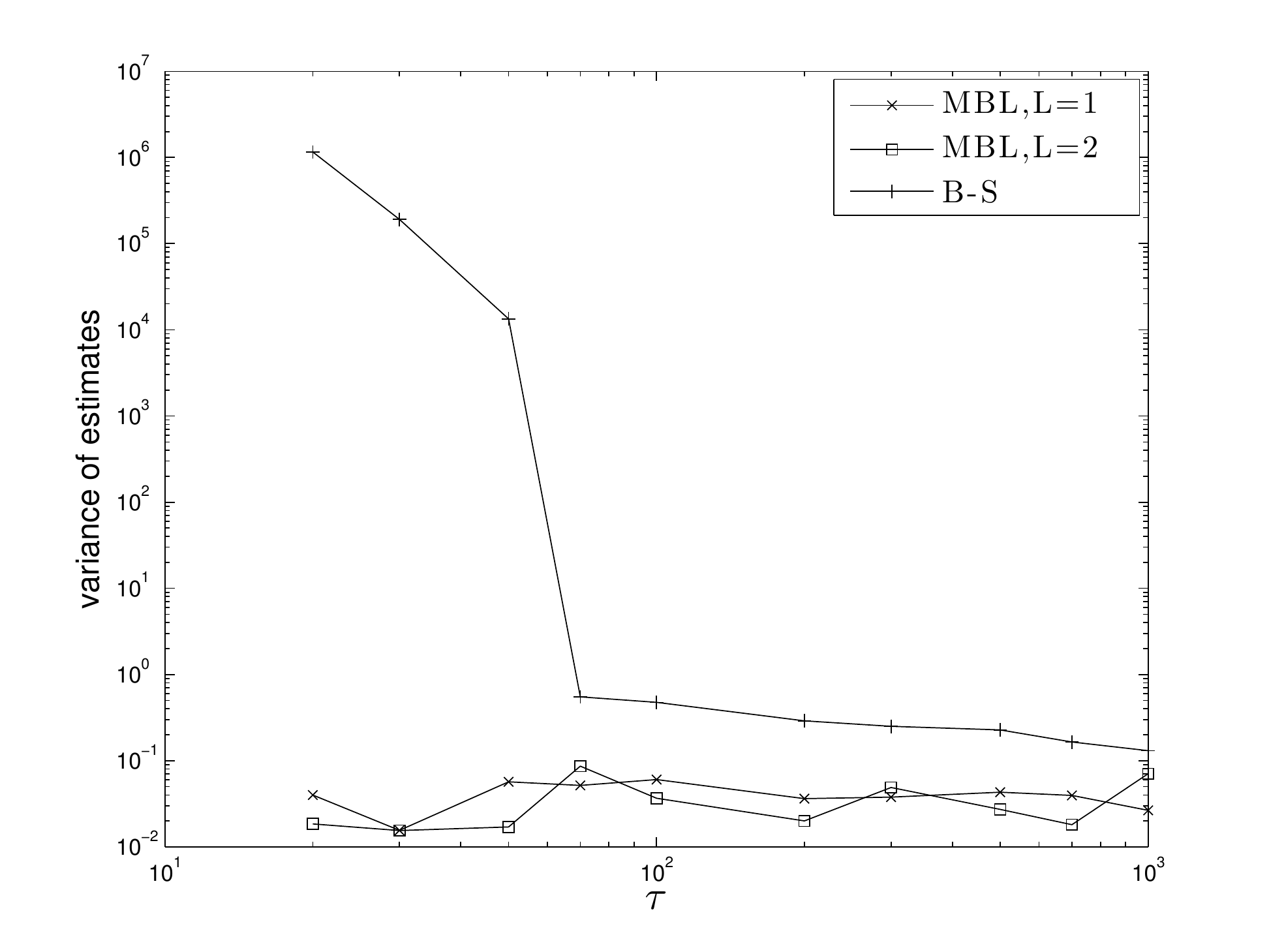}
\caption{(Linear 2D model with full observation). {\bf Left}: the mean relative root-mean-square error (MRrmse) of the estimates from Berry-Sauer's method (B-S) with $\tau=2000$, Belanger's method with $L=1$, and modified Belanger's method with $L=1,\tau=1000$. {\bf Right}: Variance of estimates for different $\tau$: Variance of estimates of modified Belanger's method (MBL) with lags $L=1, 2$ compared to the variance of estimates of Berry-Sauer's method over the last 5000 time steps.}
\label{figLin2D_2obs}
\end{figure}

In Figure~\ref{figLin2D_2obs} (left panel) we show the mean relative root-mean-square error of the estimates of diagonal components of $Q$ and $R$ as functions of time, from Berry-Sauer's method, the original Belanger's method, and the modified Belanger's method, respectively. The mean relative root-mean-square error is defined by
\BEA
MRrmse(t)=\frac{1}{n+m}(\sum_{s=1}^{n}\frac{|\widetilde{Q}_{s,s,t}-Q_{s,s}|}{Q_{s,s}}+\sum_{s=1}^{m}\frac{|\widetilde{R}_{s,s,t}-R_{s,s}|}{R_{s,s}})\label{MRrmse} \label{mrrmse}
\EEA,
where $\widetilde{Q}_{s,s,t}$ are $\widetilde{R}_{s,s,t}$ are the $s$th diagonal element of the estimates of $Q$ and $R$ at time $t$, respectively.
 In the estimation with Berry-Sauer's method, we use $\tau=2000$ so that the variation of the estimation is small enough and we can see the estimates converge around the true value. On the other hand, we use $\tau=1000$ in the estimation with the modified Belanger's method. For comparison purpose, we also use $L=1$ for  both the original and modified Belanger's method. Notice that while all three methods work reasonably well qualitatively, the errors of the modified and the original Belanger's methods are on the same order and smaller than that of the Berry-Sauer's method. In fact, the results are qualitatively similar when the number of lags are increased (results are not shown).

Now we check the sensitivity of both Berry-Sauer's method and the modified Belanger's method to the choice of the nuisance parameter, $\tau$. In the right panel of Figure~\ref{figLin2D_2obs}, we show the variances of the estimates of $Q$ and $R$ (based on averaging over the last 5000 time step estimates) as functions of $\tau$. From this numerical test, it is clear that the modified Belanger's method is not only more robust to the choice of $\tau$, it also produces estimates with smaller variances, either using innovations up to lag $L=1$ or $L=2$.

\begin{figure}[H]
\centering
\includegraphics[width=0.7\textwidth]{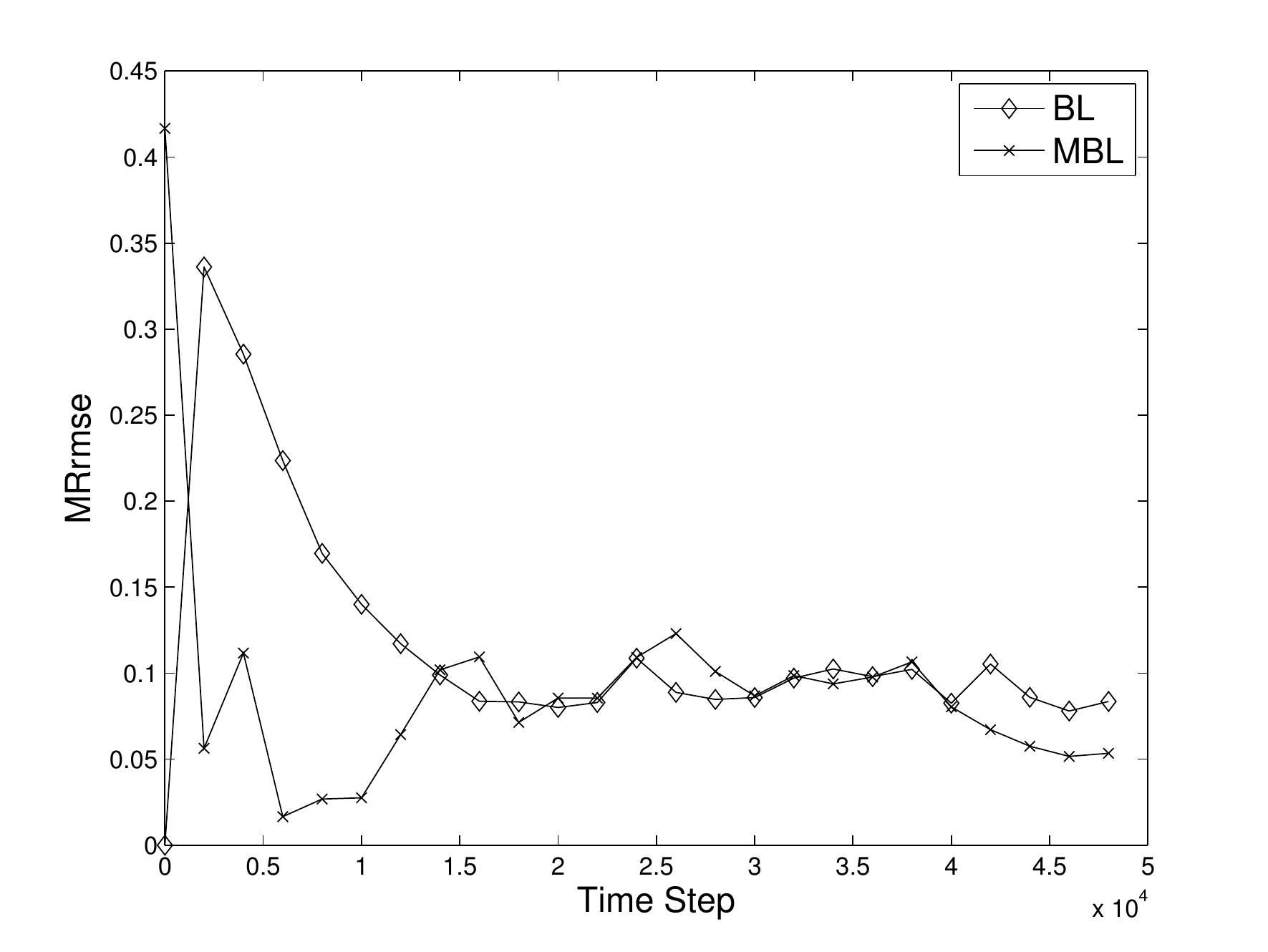}
\caption{(Linear 2D model with partial observation): The mean relative root-mean-square error of the estimates from Belanger's method with $L=4$ and modified Belanger's method with $L=4,\tau=1000$.}
\label{figLin2D_1obs}
\end{figure}
\begin{figure}[H]
\includegraphics[width=0.5\textwidth]{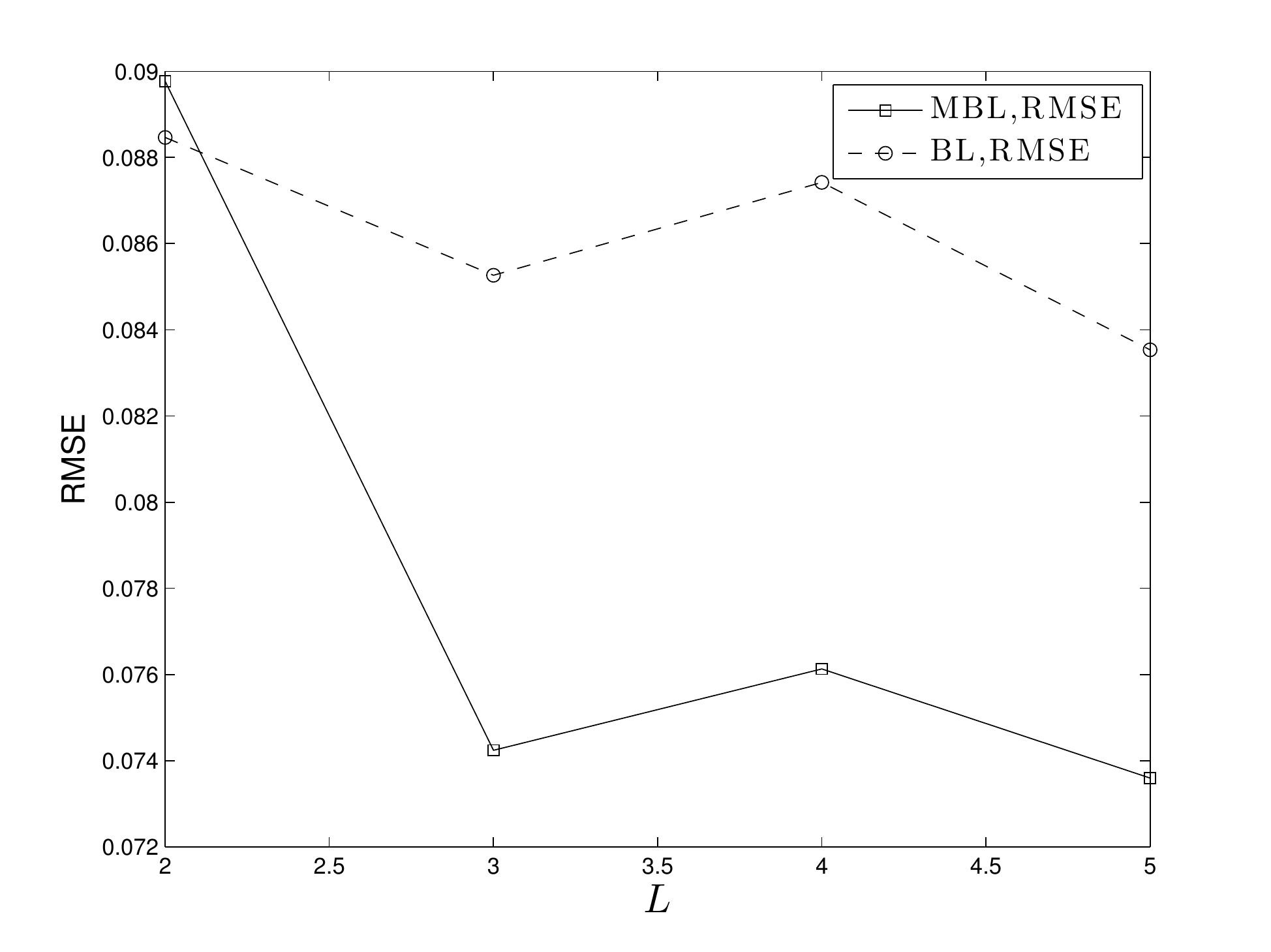}
\includegraphics[width=0.5\textwidth]{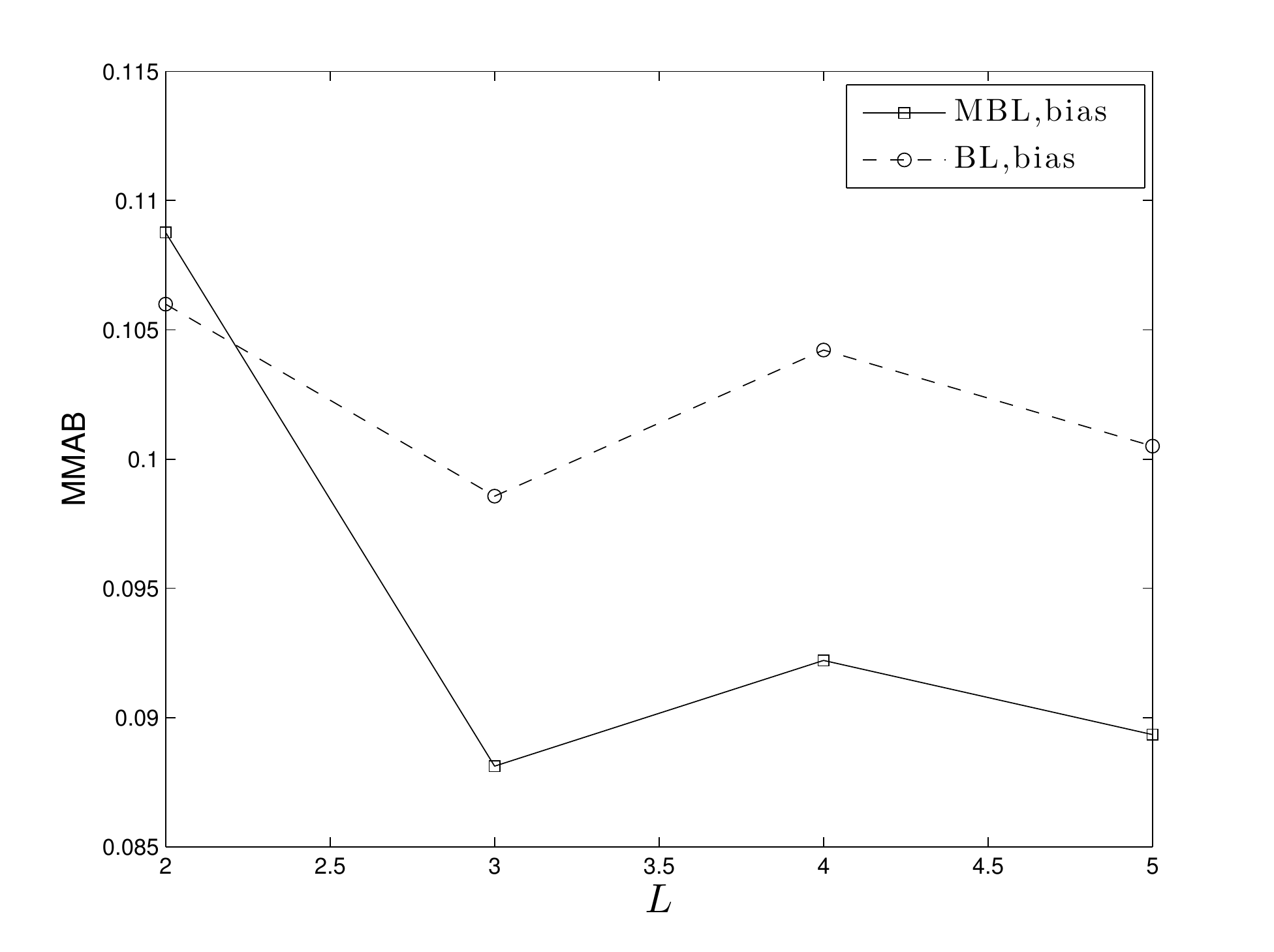}
\caption{(Linear 2D model with partial observation): Modified Belanger's method (MBL) with $\tau=1000$ and the original Belanger's method (BL) under linear 2D model with partial observations. Root-mean-square error (left) and mean of maximum absolute bias (right) of $QR$ estimates as functions of $L$. }
\label{figLin2D_1obsv2}
\end{figure}

Figure~\ref{figLin2D_1obs}  shows the comparison results of the original Belanger's method and the modified Belanger's method in the context of partial observations. Here, we show results with $L=4$ for both methods and $\tau=1000$ for modified Belanger's method. Note that the results for both methods are comparable in this experiment. In this case, Berry-Sauer's estimation method does not work (results are not shown) since the regression problem in \eqref{solveQ} is under determined. Furthermore, this regression problem can become ill-posed for some choice of $\Gamma$; e.g., when $\Gamma=\mathcal{I}$, the regression coefficient for $\alpha_2$ in \eqref{solveQ}, that is, $H_{j+1}F_j\Gamma Q_2\Gamma^{\top}H_j^\top$ is zero and therefore the second diagonal component of $Q$ is unobservable. Indeed, we tested both Belanger's scheme with various lags and found that the estimates are only accurate for $L>1$. Figure~\ref{figLin2D_1obsv2} shows the root-mean-square error (RMSE) and the mean of maximum absolute bias (MMAB) of these two versions of Belanger's method. More precisely, we perform 200 independent  simulations, each of these assimilates 50,000 steps. These two measuring skills are computed by averaging over these 200 simulations, accounting estimates from the last 5000 steps of 50,000 assimilation steps:
\BEA
\text{RMSE}&=&\sqrt{\frac{1}{200}\sum_{i=1}^{200}\frac{1}{5000}\sum_{k=45001}^{50000}\|Q_k-Q\|^2+\|R_k-R\|^2}\nonumber\\
\text{MMAB}&=&\frac{1}{200}\sum_{i=1}^{200}\max\Big\{|(\frac{1}{5000}\sum_{k=45001}^{50000}Q_k)-Q|,|(\frac{1}{5000}\sum_{k=45001}^{50000}R_k)-R|\Big\}\nonumber,
\EEA
where the maximum is over the entries of these matrices.
Notice that the estimates of both methods are accurate according to these measures and the modified Belanger's scheme produces slightly more improved estimates as $L$ increases.

\subsection{A low-dimensional nonlinear example}
In this section, we consider filtering problem of a low-dimensional system of SDE's, 
\BEA
\frac{dx}{dt}= Mx+B(x,x) -Dx+\Gamma \dot{W_t},\label{eq18}
\EEA
where $x=(u,v_1,v_2)^{\top}$, $B(x,x)=(0, auv_2,-auv_1)^{\top}$, $D=diag(0, d_1, d_2)$, $M=\begin{pmatrix}0 & \omega & 0\\-2\omega & 0 & -\theta\\ 0 & \theta & 0\end{pmatrix}$, $\Gamma=\begin{pmatrix}0 & 0\\\sigma_1 & 0\\0 & \sigma_2\end{pmatrix}$ and $W_t$ is the Wiener process with $\mathbb{E}(W_tW_t^{\top})=Qt$ where $Q=I_2$. The system in \eqref{eq18}, which was called the zeroth level memory model in \cite{HarlimMahdiMajda2014}, was derived to model the interaction of a zonal jet $u$ with topographic waves, $v_j$. In particular, the terms $Mx$ and $B(x,x)$ in \eqref{eq18} are special solutions of the topographic stress model on a two-dimensional periodic domain, truncated up to a total horizontal wavenumber-one \cite{cf:87,mw:06}, while the linear damping and noises, $-Dx+\Gamma \dot{W_t}$, terms are added to represent the interactions with the unresolved higher-order Rossby modes. The triad model in \eqref{eq18} has some interesting properties, namely, it is geometrically ergodic provided that $d_1, d_2>0$, the system has a Gaussian invariant measure, and the systems equilibrium statistics are equivalent to those when $\omega$ is replaced by $-\omega$ (see \cite{HarlimMahdiMajda2014} for details).

In the numerical tests below we set $a=1$, $\omega=\frac{3}{4}$, $\theta=1$, $d_1=d_2= \sigma_1^2=\sigma_2^2 = \frac{1}{2}$ as in \cite{HarlimMahdiMajda2014}. We choose integration time step and observation time step to be $\delta t=\Delta t=0.1$. An explicit Euler's scheme is implemented to integrate  \eqref{eq18}. Here, we implement the parameter estimation method with ETKF (see Algorithm~\ref{ETKF_k} for details) with an ensemble of size $N_e=16$. In the context of full observations,  $H_j=\mathcal{I}_3$ is the identity matrix, and $R=0.257\mathcal{I}_3\Delta t$ is a diagonal matrix (the observation noise amplitude roughly corresponds to 10\% of the variance of $u$). On the other hand, in the context of partial observations $H_j=(1,0,0)$ and $R=0.257\Delta t$ is a scalar. Since $Q$ and $R$ are diagonal in both situations, we naturally parameterize them by their diagonal components, exactly as in \eqref{param}.

\begin{figure}
\includegraphics[width=0.5\textwidth]{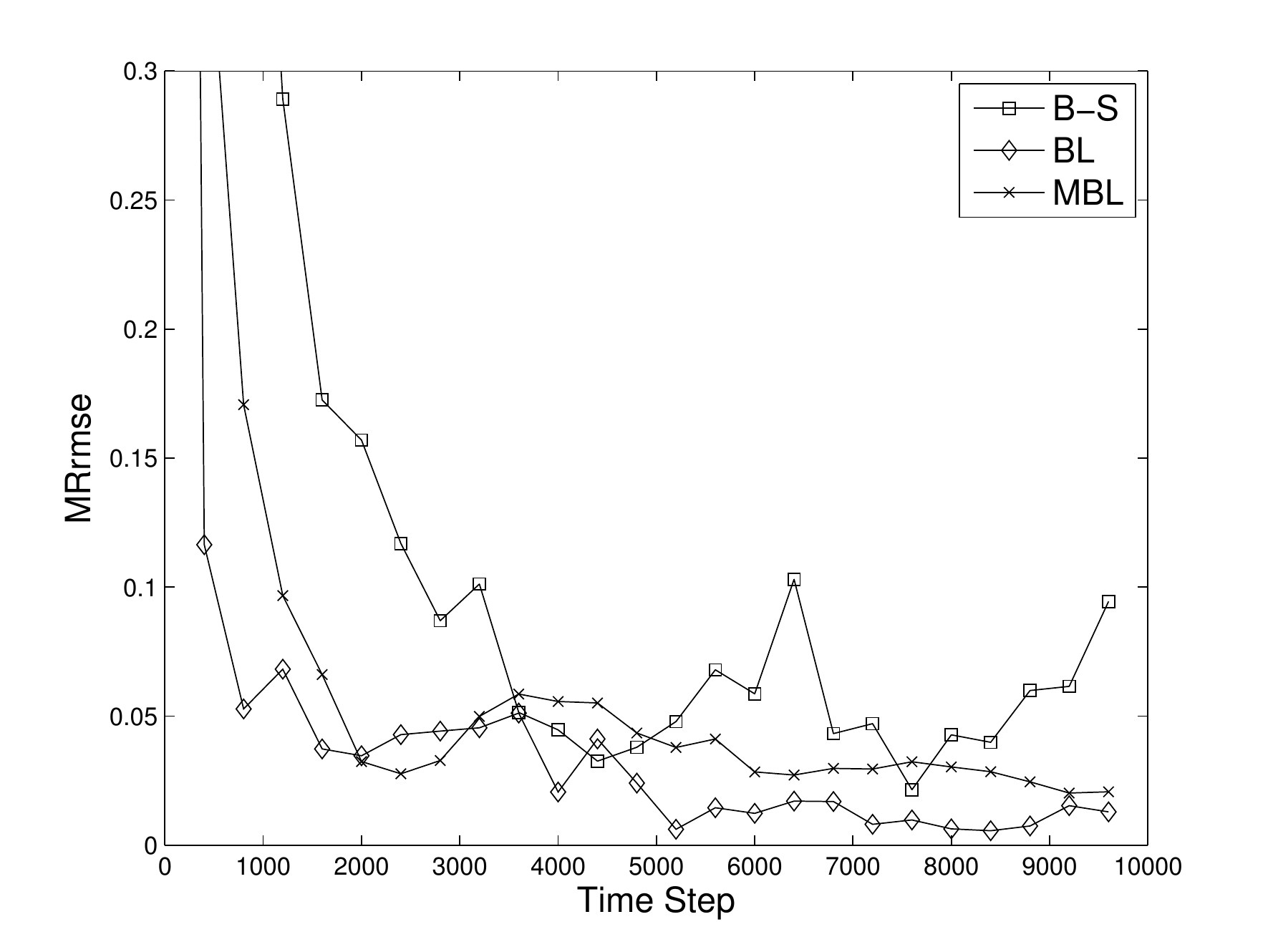}
\includegraphics[width=0.5\textwidth]{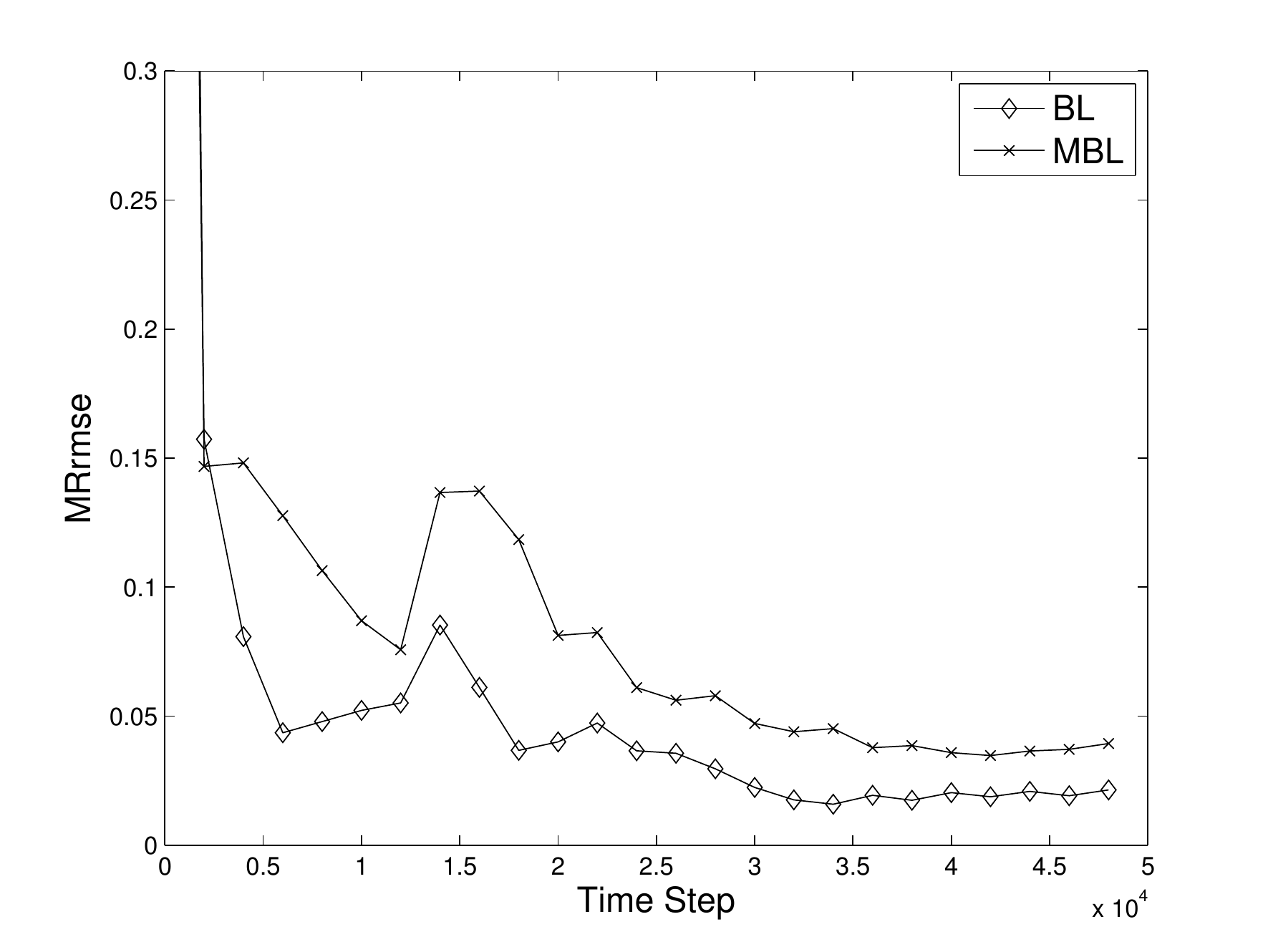}
\caption{{\bf Left:} (triad model with full observation) the mean relative root-mean-square error of the estimates from Berry-Sauer's method with $\tau=1000$, Belanger's method with $L=1$, and Modified Belanger's method  with $L=1,\tau=100$. {\bf Right:} (triad model with partial observations)  the mean relative root-mean-square error of estimates from Belanger's method with $L=8$ and Modified Belanger's method with $L=8,\tau=1000$. }
\label{figTriad}
\end{figure}

In Figure~\ref{figTriad} (left panel), we compare the performance in the context of full observations. Notice that the relative errors of all three methods are comparable. For this experiment, we set $\tau=100$ for the modified Belanger's method and $\tau=1000$ for Berry-Sauer's method so the estimates have small variances. For a fair comparison, we also set $L=1$ for both the original and modified Belanger's methods. While the relative errors in both versions of Belanger's methods decay faster than that of Berry-Sauer's method, notice that the relative error of the original Belanger's method decays faster than that of the modified Belanger's method. We speculate that this is because the cost function in \eqref{belanger} is defined with a more appropriate norm, as mentioned in Section~2.1.2.
 
In the context of partial observations, only the first component of the state variables is observed. Again we can not solve Berry-Sauer's least square problem in \eqref{solveQ} since it is underdetermined. The right panel of Figure~\ref{figTriad}  shows the estimates from the two versions of Belanger's method with lag $L=8$ and $\tau=1000$. Notice that the relative errors of both methods decay slower when only partial observations are available. 

\subsection{A stochastically forced 40-dimensional Lorenz-96 example}

In this section, we compare the proposed modified Belanger's method only to Berry-Sauer's method on a relatively higher-dimensional problems, the 40-dimensional Lorenz-96 model \cite{Lorenz1996}, which has been routinely used to validate data assimilation schemes. We ignore showing numerical results with the original Belanger's method since the computational cost significantly increases as the number of observations, $m$, becomes large. In our experiment, we follow \cite{BerrySauer2013}, forcing the Lorenz-96 model \cite{Lorenz1996} with Gaussian white noises. In particular, the system of the stochastically forced Lorenz-96 model is given by:
\BEA
\frac{dx_i}{dt}=-x_{i-2}x_{i-1}+x_{i-1}x_{i+1}-x_i+F+\Gamma_{i}\dot{W_t},\label{L96}
\EEA
where $\Gamma_i$ denotes the $i$th row of the matrix $\Gamma$, which we will set as an identity matrix in our numerical experiments below, and the Wiener process $W_t$ has covariance $\mathbb{E}(W_tW_t^{\top})=\hat{Q}t$. We use the standard setting for this model with $F=8$ and integration step size, $\delta t=0.05$. The deterministic part of \eqref{L96} is integrated with a fourth-order Runge-Kutta method while the stochastic part is integrated with the standard Euler's scheme.  In our numerical experiments, we consider sparse, 20 equally spaced, observations of the state variables for every $N$ integration time steps, i.e. $\Delta t=N\delta t$.  We generate the true covariance matrix $\hat{Q}$ randomly in the way that the eigenvalue of $\hat{Q}$ distribute uniformly between $0.1$ and $1$. The observation error covariance $R$ is generated in the same way, but multiplied by a scalar so that we can control the ratio of the error covariances, $\text{tr}(R)/\text{tr}(Q)$, where $Q=\hat{Q}\delta t$ is what we will estimate. In particular, when $\text{tr}(Q)=\text{tr}(R)$ the magnitude of $R$ is about $0.025$ which is  about of the same magnitude as that used in \cite{ZhangZhangHansen2009}. 
For this nonlinear example, we will apply the modified Belanger's scheme with the ETKF algorithm~\eqref{ETKF_k} with an ensemble of size $N_e=50$ and a total number of $50000$ assimilation steps.   We parameterize $Q$ and $R$ as in \eqref{param}, dividing matrices the matrix $Q$ into $4\times 4$ blocks: 
\BEA
Q=\begin{pmatrix}\mathcal{B}_{1,1} & \mathcal{B}_{1,2} &... & \mathcal{B}_{1,10}\\\mathcal{B}_{2,1} & \mathcal{B}_{2,2} & ... & \mathcal{B}_{2,10}\\\vdots\\ \mathcal{B}_{10,1} & \mathcal{B}_{10,2} & ... & \mathcal{B}_{10,10}\end{pmatrix}, \text{and let } Q_s=\begin{pmatrix}0 & ... & ... & ...& ...& 0\\\vdots \\0 & ... &... & \mathcal{B}_{s_i,s_j} & ... & 0\\\vdots \\0 & ...& \mathcal{B}_{s_j,s_i} & ... & 0 & 0\\\vdots\\ 0 & ... & ... & ... & ... &0\end{pmatrix},\nonumber
\EEA 
where $s_i$ and $s_j$ are proper subindexes corresponding to $s$. Hence there are $N_Q=55$ parameters for $Q$. On the other hand, we do not parameterize $R$, resulting to $N_R=210$ parameters for $R$ by symmetry. When we estimate $\alpha_s$ and $\beta_s$,  $\tau=10000$ is used for Berry-Sauer's method and $\tau=1000$ and $L=1$ or $L=3$ are used for the modified Belanger's method.

\begin{figure}
\includegraphics[width=0.9\textwidth]{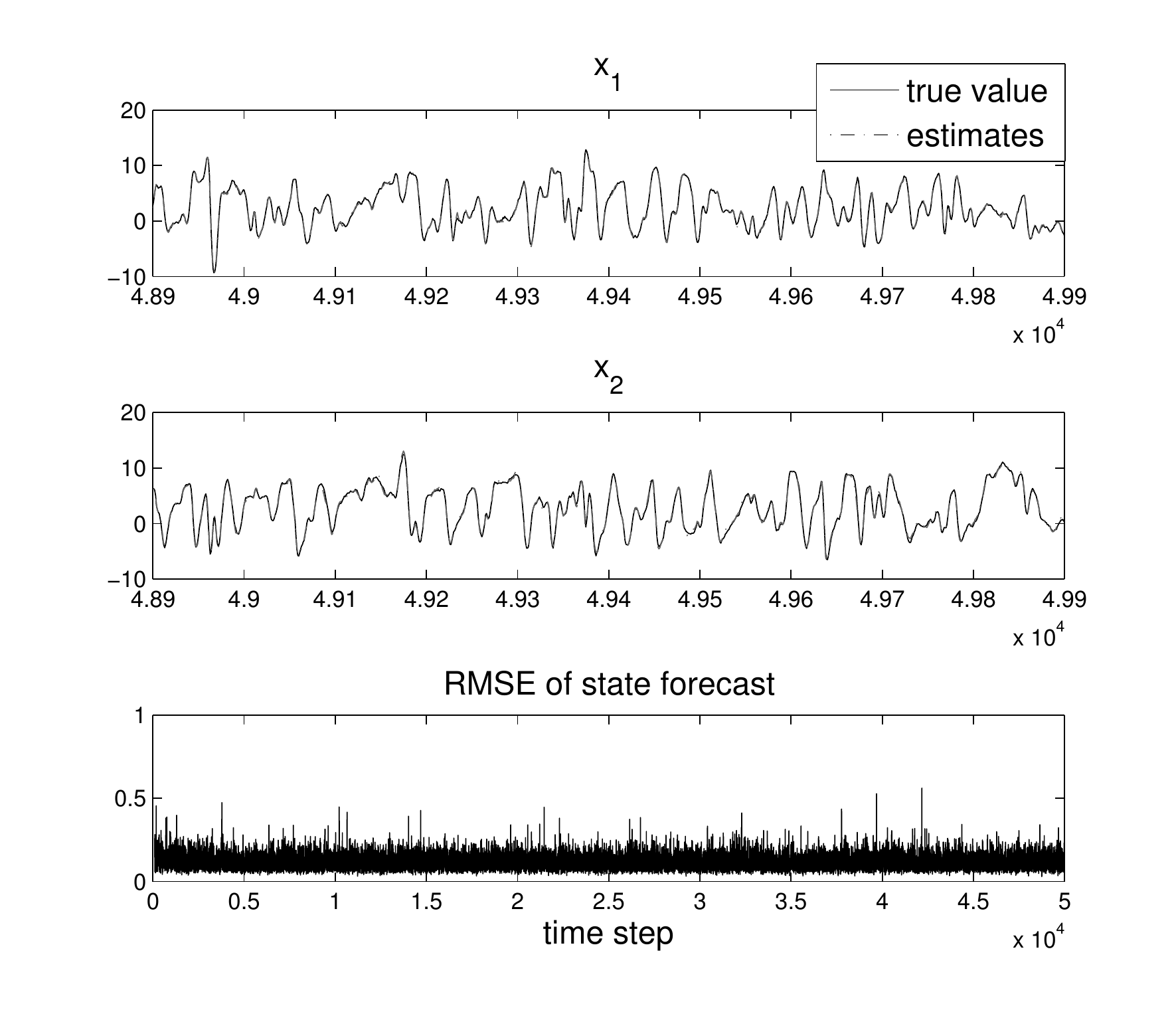}
\caption{Modified Belanger's scheme: Estimates of $x_1$ (observed variable) and $x_2$ (unobserved variable) from the 48900$th$ time step to the 49900$th$ time step, and RMSE of the analysis state.}
\label{figL96_20obs_MBL_x}
\end{figure}

Figure~\ref{figL96_20obs_MBL_x} shows that reasonably accurate estimates of the first two components  $x_1$ (one of the observed components) and $x_2$ (one of the unobserved components) of the state variable can be obtained by ensemble transform Kalman filter incorporated with modified Belanger's method with $N=1$, $L=3$, and $\text{tr}(R)/\text{tr}(Q)=1$. We also include the RMSE of the analysis state.

%

\begin{figure}[H]
\includegraphics[width=0.5\textwidth]{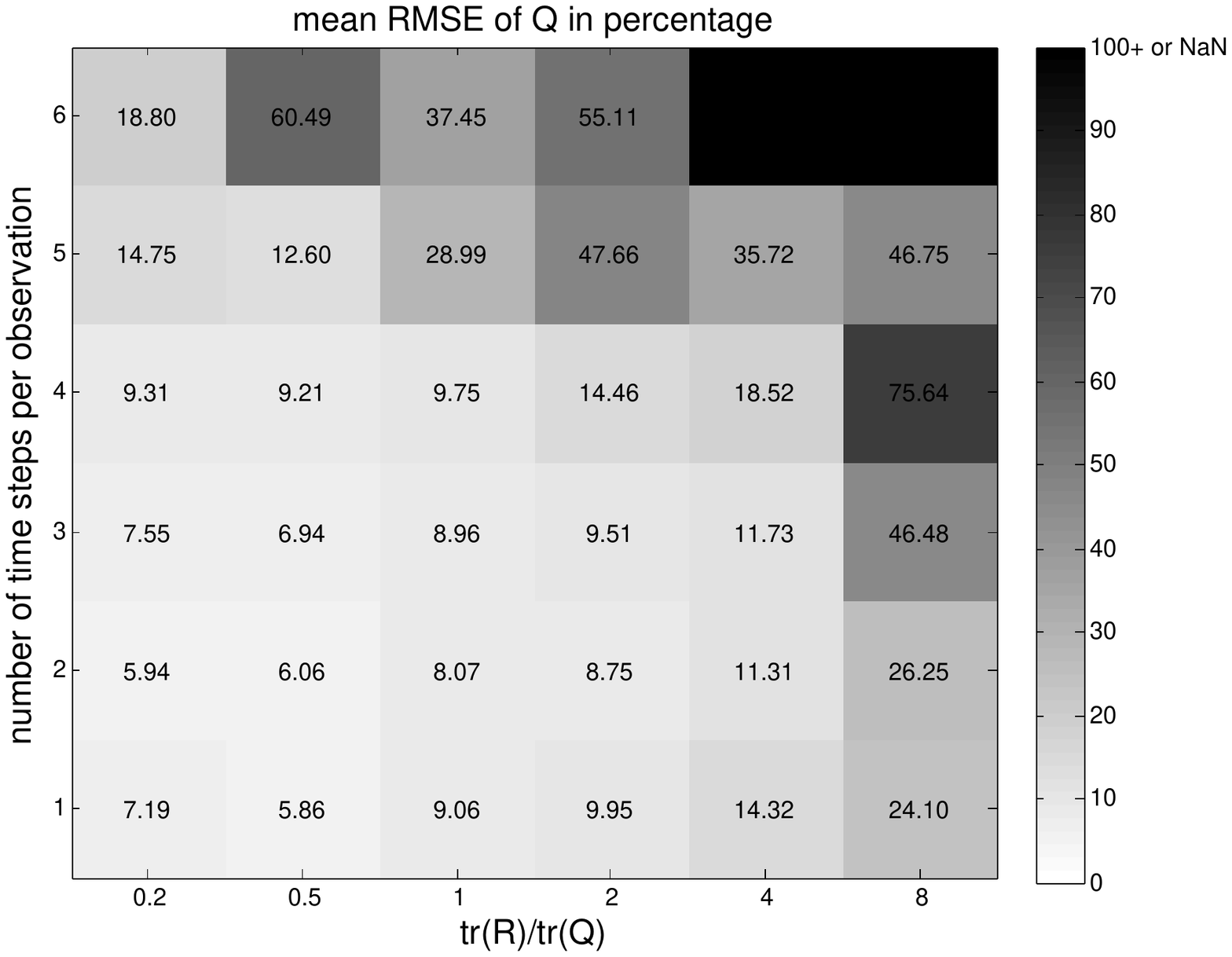}
\includegraphics[width=0.5\textwidth]{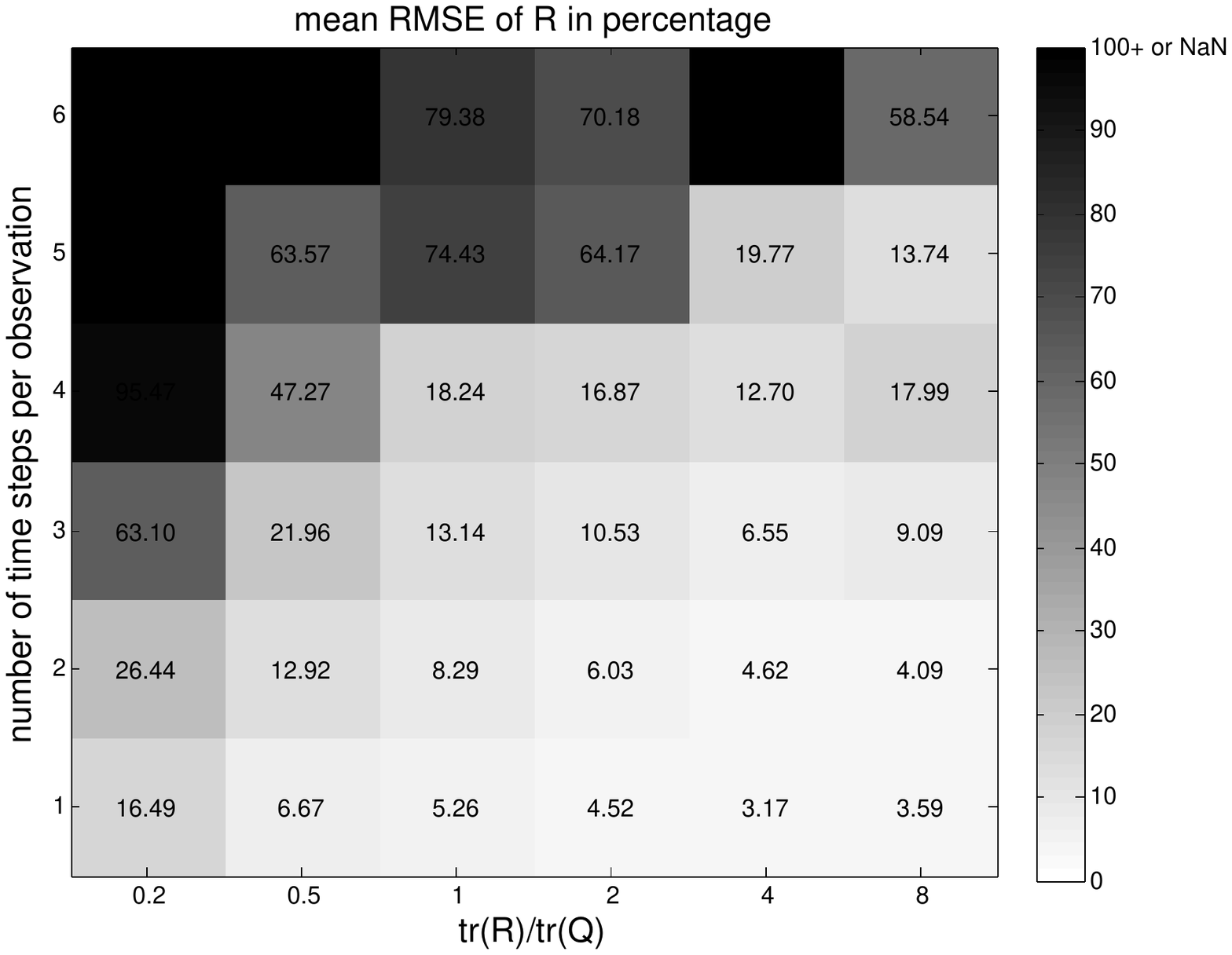}
\caption{Modified Belanger's method with $L=1,\tau=1000$ and Lorenz-96 model: the average of the root-mean-square error (in percentage) of the estimates of $Q$ (left) and $R$ (right) over the last 25000 time steps.}
\label{figL96_L1_20obs_MBL_QR}
\end{figure}

\begin{figure}[H]
\includegraphics[width=0.5\textwidth]{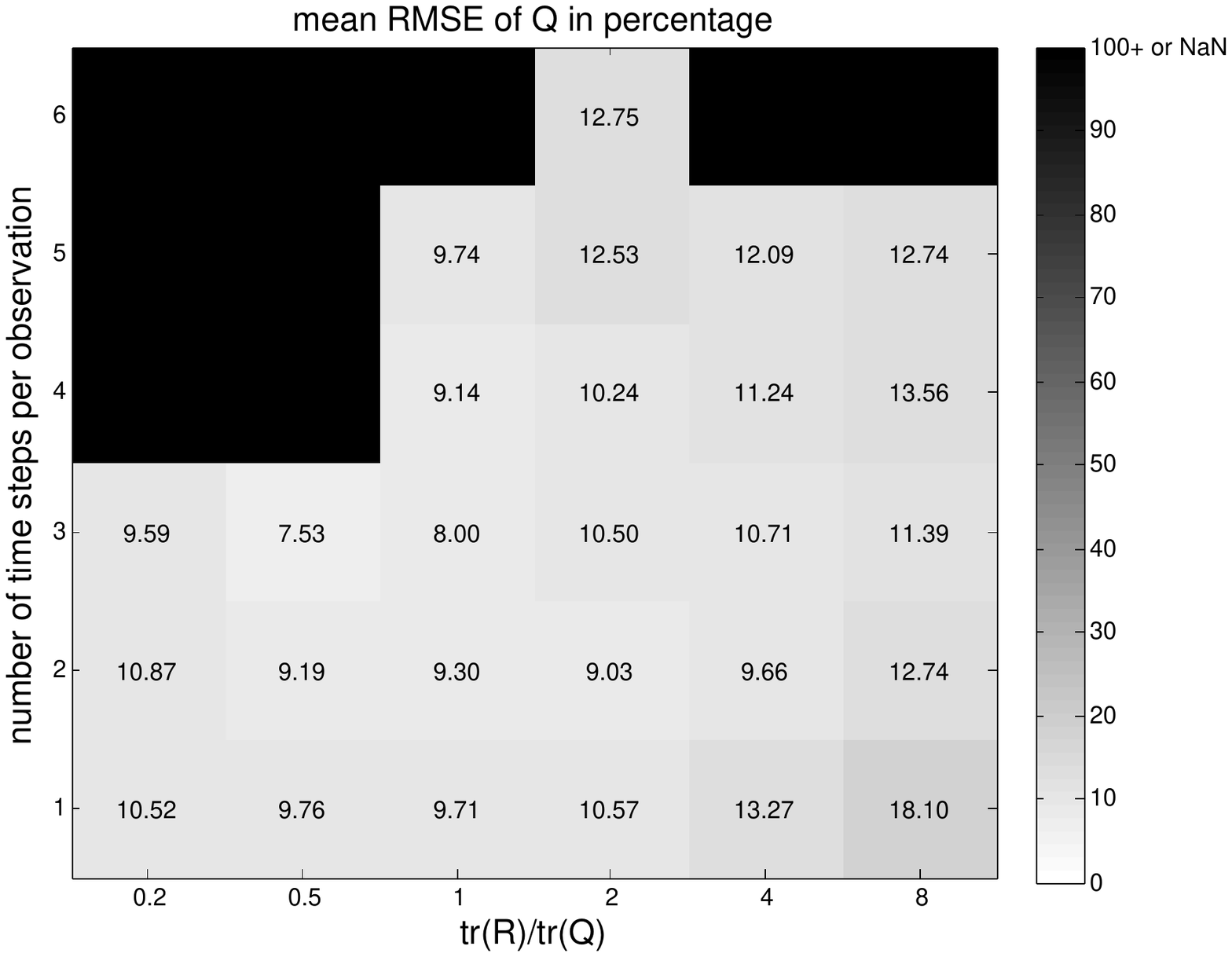}
\includegraphics[width=0.5\textwidth]{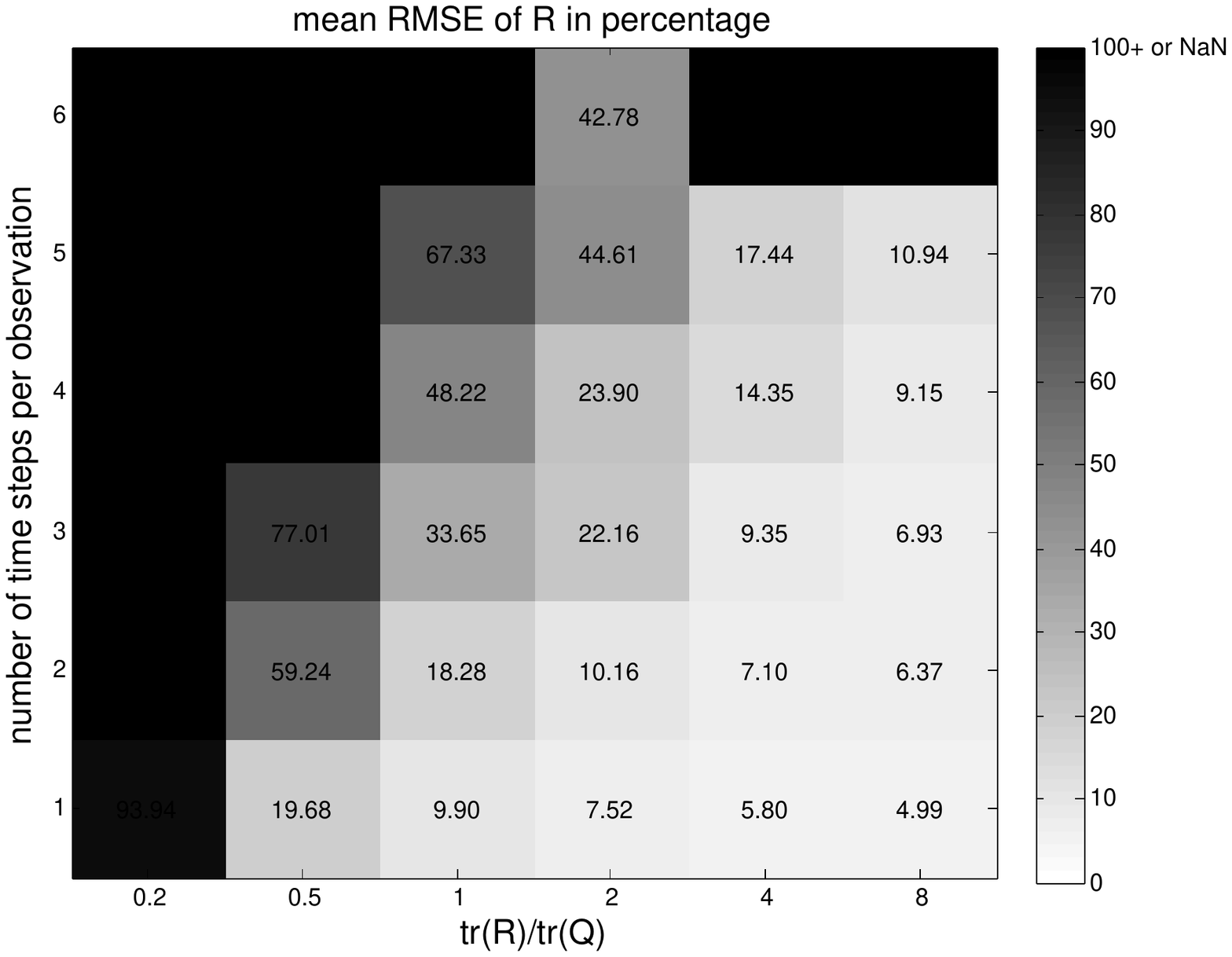}
\caption{Berry-Sauer's method with $\tau=10000$ and Lorenz-96 model: the average of the root-mean-square error (in percentage) of the estimates of $Q$ (left) and $R$ (right) over the last 25000 time steps. }
\label{figL96_20obs_BS_QR}
\end{figure}

\begin{figure}[H]
\includegraphics[width=0.5\textwidth]{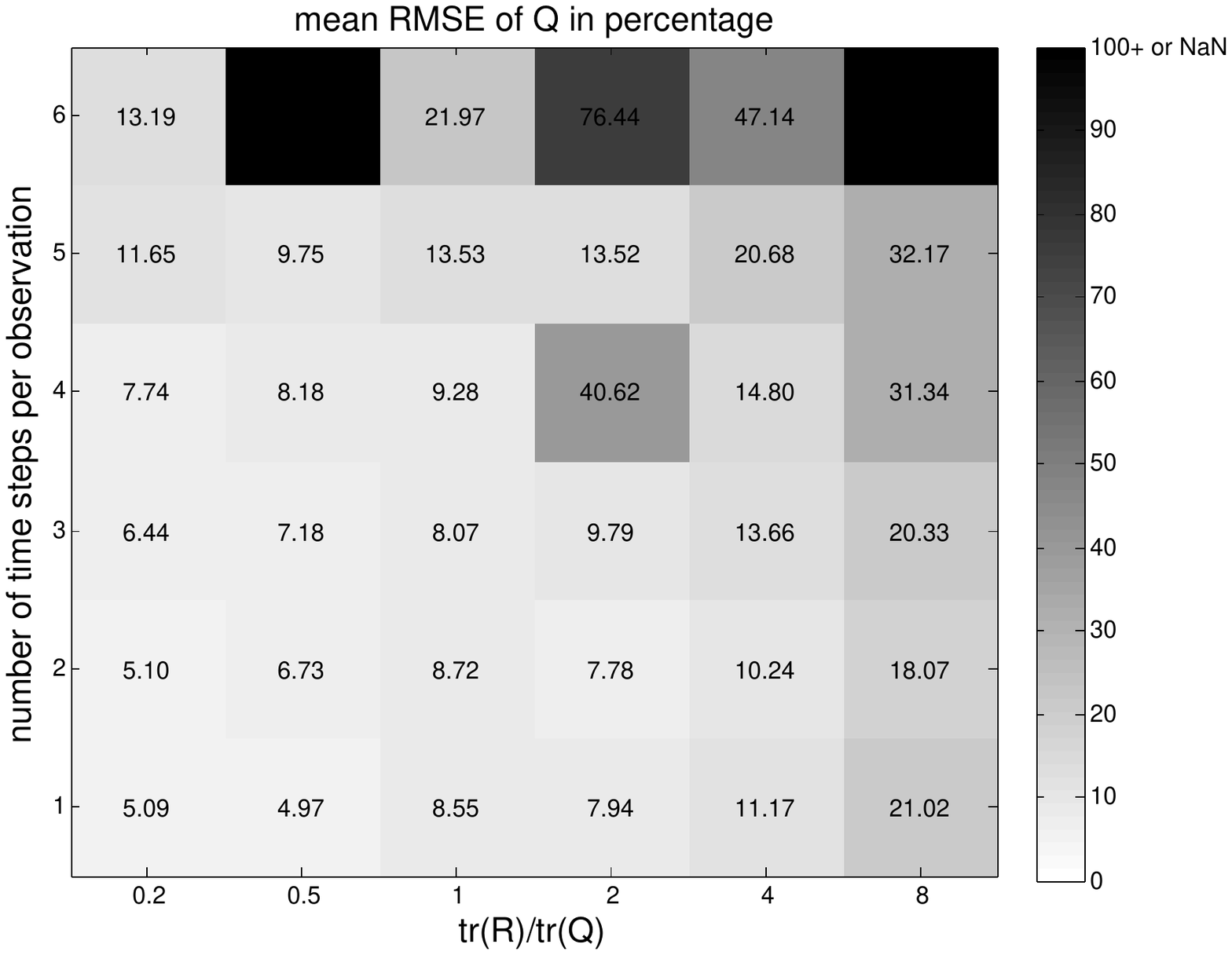}
\includegraphics[width=0.5\textwidth]{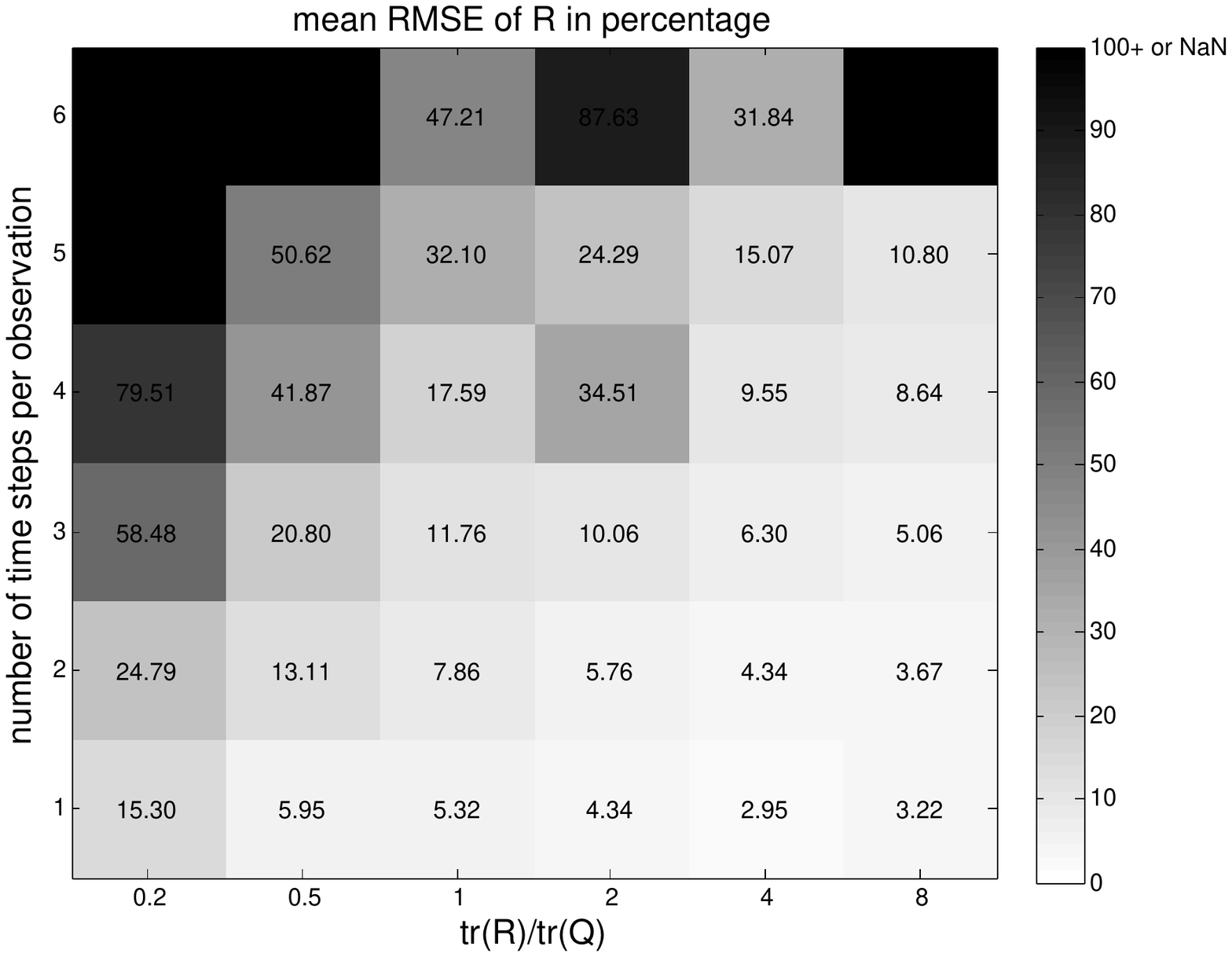}
\caption{Modified Belanger's method with $L=3,\tau=1000$ and Lorenz-96 model: the average of the root-mean-square error (in percentage) of the estimates of $Q$ and $R$ over the last 25000 time steps.}
\label{figL96_L3_20obs_MBL_QR}
\end{figure}

Next, we will discuss the filter performance for different observation step size, $N$, and covariance ratio, $\text{tr}(R)/\text{tr}(Q)$.
Figures \ref{figL96_L1_20obs_MBL_QR}-\ref{figL96_L3_20obs_MBL_QR} show the distribution of relative errors in $Q$ and $R$ based on estimates from the two methods with different $L$ (for modified Belanger's method only), respectively.  Specifically, the relative error is defined by the average over the last 25000 time steps of the following error in percentage:
\BEA
\text{error in percentage}=\frac{\|\text{estimates}-\text{truth}\|}{\|\text{truth}\|}\times 100,
\EEA
where the \lq\lq estimates" and \lq\lq truth" refer to the estimation of $Q$ (or $R$, respectively) at a specific time step and the true $Q$ (or $R$, respectively) matrix and $\|\cdot\|$ is  the Frobenius norm. Figure \ref{figL96_L1_20obs_MBL_QR} shows the 
relative errors from the modified Belanger's method with $L=1$. We notice that the estimations of $Q$ and $R$ have different sensitive regions. More generally, for fixed observation time step $\Delta t=N\delta t$, the estimation of $Q$ is less accurate when the magnitude of $R$ is larger, and vice versa. This result can be understood as follows. Recall that the observational matrices $\mathcal{H}^{(Q)}$ (or $\mathcal{H}^{(R)}$) in \eqref{eq6} (or cf. \eqref{HQjls}, \eqref{HRjls}) for the parameter space are linear functions of the basis of $Q$ (or $R$). Since we parameterize $Q$ (or $R$) block-wise by the true value of the sub-blocks of  $Q$ (or $R$),  larger magnitude of $Q$ (or $R$) causes $\mathcal{H}^{(Q)}$ (or $\mathcal{H}^{(R)}$)  to have larger entries. Since we solve equation \eqref{newmethod} directly by a least-square method, we naturally expect more accurate estimates for the parameters with larger coefficients. 

We should note that it is possible to use a more appropriate matrix norm  in \eqref{newmethod} that is similar to the one proposed in the original Belanger's method, so that we can avoid the estimation of $Q$ and $R$ to have different precision when their magnitudes are significantly different. However this may also raise the computational cost to the same order of that of the original Belanger's method, which is what we are trying to avoid. According to the analysis above, it is reasonable to deduce that good estimates of $Q$ can be achieved by the modified Belanger's method when the magnitude of $R$ is small relative to $Q$, and vice versa. 

On the other hand, we also found that the estimates of the diagonal elements of $Q$ and $R$ are more accurate in most scenarios. We should note that we have also tested on much larger observation noise error with diagonal component $1$ and our method still converges (results are not reported). In fact, this method can be applied to much broader choices of $N$ and $\text{tr}(R)/\text{tr}(Q)$ if the true $Q$ and $R$ are diagonal, which motivates one to find better choices of $Q_s$ and $R_s$ to reduce the sensitivity of the estimates when $Q$ and $R$ have different magnitudes. Furthermore, better choices of $Q_s$ and $R_s$ can reduce the number of unknown parameters, $\alpha_s, \beta_s$, which is particularly important for large dimensional problems. In our numerical experiments, none of the methods can work when the number of parameters are larger than the number of product of innovation processes that are being used for the covariance estimation. 

In Figure~\ref{figL96_20obs_BS_QR}, we show the corresponding relative errors of Berry-Sauer's method. We see that the estimation of $Q$ and $R$ by the method of Berry-Sauer's is less stable when $\text{tr}(R)/\text{tr}(Q)$ is small. However the estimation of $Q$ of their method is more accurate when $\text{tr}(R)/\text{tr}(Q)$ is large. Notice also that the resulting estimates of $Q$ and $R$ do not have the opposite sensitivity when the ratio of $\text{tr}(R)/\text{tr}(Q)$ are different than 1 as encountered in the modified Belanger's scheme. We suspect that this may be because they estimate $R$ in \eqref{solveR} and $Q$ in \eqref{solveQ} separately.

Figure~\ref{figL96_L3_20obs_MBL_QR} shows the relative errors of the modified Belanger's method with larger lags, $L=3$. We see that in many cases the errors indeed can be reduced by incorporating more lagged information. For example, the relative errors of the estimates of $R$ can be reduced from $74.49\%$ to $32.10\%$ when $L$ is increased from 1 to 3, for $N=5$ and $\text{tr}(R)/\text{tr}(Q)=1$. However, we should also mention that increasing lags do not seem to improve the estimates in some regimes with large $N$ (e.g., $N=6$ and $\text{tr}(R)/\text{tr}(Q)=0.2$).

\begin{figure}[H]
\includegraphics[width=0.7\textwidth]{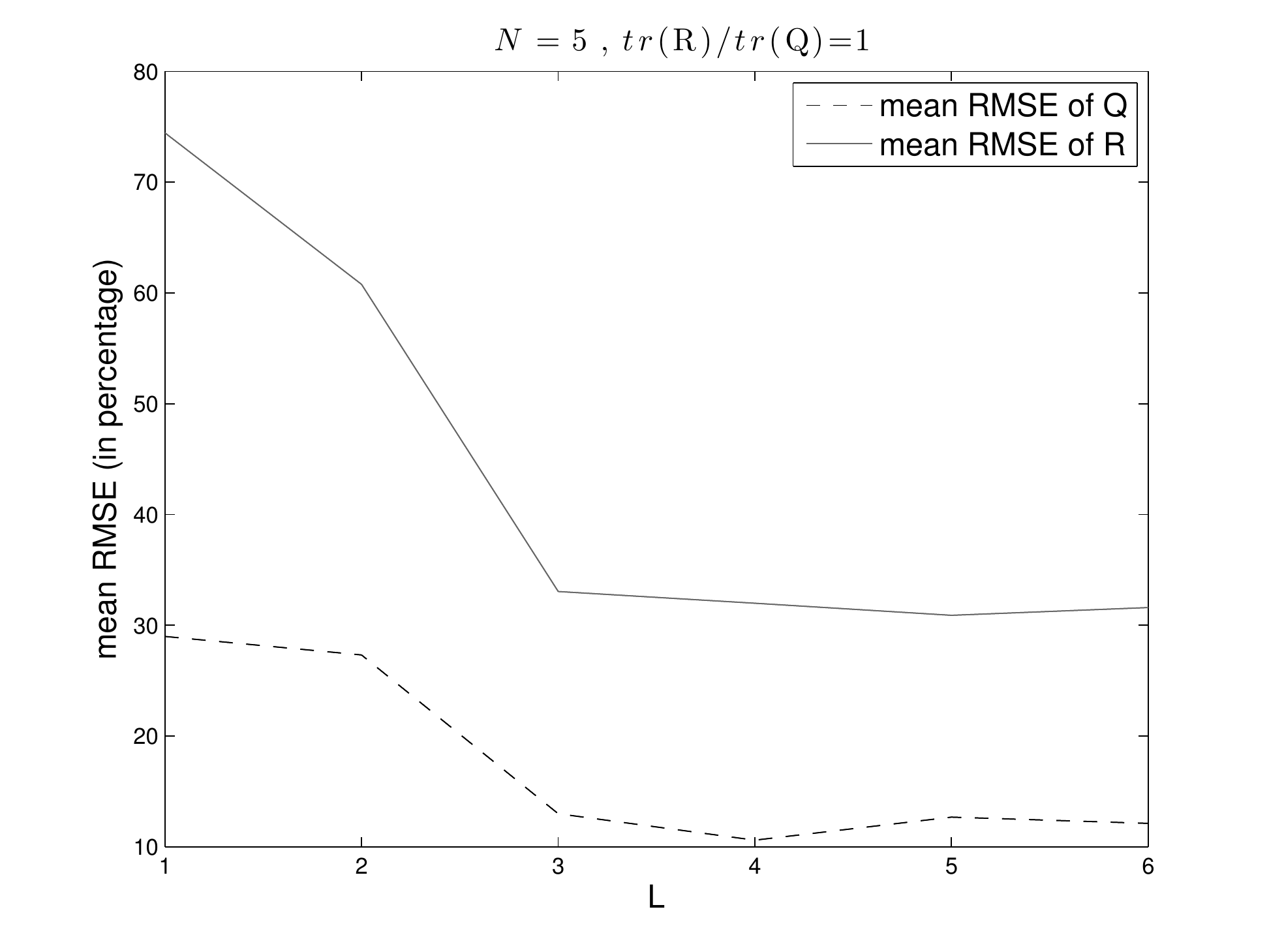}
\caption{Modified Belanger's method with $L=1,2,3,4,5,6$ and Lorenz-96 model: the average of the root-mean-square error (in percentage) of the estimates of $Q$ and $R$ over the last 25000 time steps. }
\label{figL96_LL_20obs_MBL_QR}
\end{figure}

In Figure~\ref{figL96_LL_20obs_MBL_QR} we show the plot of the root-mean-square errors of $Q$ and $R$ as functions of $L$, for the regime $N=5$ and $\text{tr}(R)/\text{tr}(Q)=1$, a case for which the method of Berry-Sauer's does not work well as can be seen in Figure~\ref{figL96_20obs_BS_QR}. Again, we find that the error of estimates can be reduced by increasing $L$. But the improvement is small starting from $L=3$.  We suspect that this is because the entries, that corresponds to lags greater than or equal to 3, in the coefficient matrix $\mathcal{H}^{(Q)}$ and $\mathcal{H^{(R)}}$ are small, as we numerically verified. Hence adding more lags beyond $L=3$ does not increase too much solvability of equation \eqref{newmethod}.  

\subsection{Applications on L-96 model using LETKF}
In this section we show the numerical results of incorporating the modified Belanger's method in the local ensemble transform Kalman filter (LETKF) (\cite{ott2004,hunt:07}). A brief description of how to apply the modified Belanger's method in the context of LETKF was described in Section~2.2. For comparison purpose, we also show results of implementing the original Belanger's method with LETKF in a similar way. Here we do not show the results with Berry-Sauer's method since they already did it in their paper \cite{BerrySauer2013}. We should also note that the incorporation of these methods on other types of localization schemes  (\cite{GaspariCohn1999,JAnderson2011,LeiAnderson2014,BishopHodyss2007,ZhenZhang2014,Flowerdew2015}, etc.) can be nontrivial and beyond the scope of this paper.

In this example, we assume the true model is the original deterministic L96 model in \eqref{L96} without stochastic noise, $\Gamma_i=0$. Our goal is to test the potential of the proposed method as an adaptive additive covariance inflation method to mitigate errors due to localization and small ensemble size in addition to estimating $R$ on-the-fly. In particular, one can think of the additive inflation $Q$ as the covariance of a Wiener noise stochastic forcing in the forecast model:
\BEA
\frac{dx_i}{dt}=-x_{i-2}x_{i-1}+x_{i-1}x_{i+1}-x_i+F+\dot{W_t},\nonumber
\EEA
where we parameterize:
%
\BEA
Cov(W_{t+\Delta t}|W_{t})=:Q=\begin{pmatrix}
q_1 & q_2 & 0 & \dots & 0 & q_2 \\
q_2 & q_1 & q_2 & \dots & 0 & 0 \\
\vdots & \ddots  & \ddots & \ddots & \vdots & \vdots \\
0 & 0 & 0 & \ldots & q_1 & q_2 \\
q_2 & 0 & 0 & \dots  & q_2 & q_1
\end{pmatrix}.
\EEA 
In this experiment we assume that observations are taken at every site for every $\Delta t=0.05$, with observational error $R=\mathcal{I}$. To estimate $R$, we parameterize $R=r\mathcal{I}$, such that the truth is $r=1$. In total, we have three parameters for $Q$ and $R$; $\{q_1, q_2, r\}$. In the following numerical experiments, we set the LETKF with localization radius of 5, i.e., there are 11 sites in each local region.
We run the test for 2000 time steps for two ensemble sizes $N_e$: larger than the local state space with $N_e=20$ and smaller than the local state space with $N_e=6$. 

The temporal mean RMSE of analysis (MRMSE) and the final estimates of the parameters for $Q$ and $R$ from both the original Belanger's method (BL) and the modified Belanger's method (MBL) are shown in Table \ref{Tab1}.
\begin{table}[h]
\hspace{-2mm}
\begin{tabular}{llllllllll}
 & {\footnotesize MRMSE} & $q_1$ & $q_2$ &  r \\ \hline
 BL    & 0.28   & 0.01 & 0.01  & 1.01 \\ \hline
 MBL & 0.23   & 0.01   & 0.01 & 1.01
\end{tabular}
\hspace{6mm}
\begin{tabular}{lllllllll}
 & {\footnotesize MRMSE} & $q_1$ & $q_2$ &  r \\ \hline
 BL    & 0.77   &  0.51  & 0.04 & 0.9\\ \hline
 MBL & 0.81   & 0.71   & -0.06 & 0.59
\end{tabular}
\caption{Covariance estimates with LETKF with ensemble sizes $N_e=20$ (left) and $N_e=6$ (right).  }
\label{Tab1}
\end{table}
It can be seen that the results from both methods are comparable. Compare to the experiments with smaller ensemble size $N_e=6$, the tests with larger ensemble size $N_e=20$ produce more accurate state estimates (with smaller MRMSE) and also more accurate estimates of $R$. Here the estimates for $Q$ are small since the true model is deterministic and the ensemble size is large enough to resolve every possible direction of local error. 

In the case of $N_e=6$, the filter is confronted with a huge sampling error issue, because the ensemble size is much smaller than the dimension of the state variable of each local region. From the numerical results we observe that both methods effectively provide a multiplicative covariance inflation by underestimating $R$ and an additive covariance inflation with nontrivially larger $Q$. While the MRMSE of the state variables are relatively comparable and still less than the observation error, $1$, the weights of these inflations from the two methods are different; the original Belanger scheme provides less overall inflations compared to the modified Belanger's scheme. We suspect that these differences are due to different norms and cost functions that are minimized in the secondary filters.

\section{Summary}
We presented a modified version of Belanger's method to adaptively estimate the system noise covariance $Q$ and observation noise covariance $R$, which avoids expensive computational cost of the original Belanger's scheme. We embedded this method into the ensemble transform Kalman filter algorithm as in \cite{HarlimMahdiMajda2014} for filtering nonlinear problems. Similarities and differences between the modified Belanger's method and those in  \cite{Belanger1974,BerrySauer2013,OdelsonRajamaniRawlings2006} are discussed. Computationally, the modified Belanger's scheme applies a cheaper secondary filter; it produces estimates of $Q$ and $R$ based on running averages of linear regression solutions of a new cost function, replacing the secondary recursive Kalman filtering. Although the running average step is similar to the one proposed in \cite{BerrySauer2013}), our method is more flexible since its cost function can include product of innovation processes of more than one-lag.

From our numerical results we conclude as follows: While the original Belanger's scheme produces accurate estimates with faster convergence rate, it is computationally impractical for high-dimensional problems. On low-dimensional problems, we demonstrate that the accuracy of the estimates based on the proposed modified Belanger's method may be competitive (or even slightly better than) with those of the original Belanger's method. On the other hand, Berry-Sauer's method produces the least accurate estimates compared to the two methods even with full observations. For sparse observations, Berry-Sauer's method sometimes does not work simply because it only uses product of up to one-lag innovation processes. While both the modified Belanger's and Berry-Sauer's schemes requires one to prescribe a relaxation coefficient $\tau$, we found that the modified Belanger's scheme is not only more robust to variations of this nuisance parameter, but it also produces estimates with smaller variability.

In our numerical test with the Lorenz-96 model we showed that the modified Belanger's method has wider regimes of accurate estimation compared to Berry-Sauer's method. In particular, the modified Belanger's method is more accurate when observation noise amplitude is small relative to system noise amplitude and when observation time step is not too large. Berry-Sauer's method, on the other hand, has smaller regimes of accurate estimation but it is particularly very accurate when observations noise amplitude is large relative to system noise amplitude. Finally, we also compared the modified and original Belanger's methods using the local filter with LETKF \cite{hunt:07} and tested the potential of using these estimation schemes as an adaptive covariance inflation method in addition to an algorithm for estimating $R$. We found that their performances are comparable. Both methods significantly reduce the errors even when we use small ensemble size smaller than the dimension of the local state. This encouraging result shows the potential for high-dimensional applications.  

Persisting issues for practical implementation of both methods are appropriate choices of $Q_s$ and $R_s$ to reduce the numerical sensitivities and the number of parameters to be estimate. We also numerically found that adding more lags help improving the estimates of the new method in most regimes, however, there seems to be a limit in the improvement of the accuracy of the estimates beyond certain lags. A mathematically more careful analysis to improve our understanding of the approximation property of the modified Belanger's method is challenging and we leave this out for future studies.

Finally, we should mention that while the scheme was tested in the context of no model error, it will be interesting to see how this scheme will perform in the presence of model error.  Obviously, when model error is presence, an additional difficulty is to remove biases (or mean model error) from the estimates \cite{dds:98} and recent study suggested that accurate estimation can be obtained when one simultaneously estimates  both the mean and covariance error statistics \cite{BerryHarlim2014}, assuming that one knows how to prescribe an appropriate stochastic model as a model error estimator. We plan to test this method to address this issue in our future study.

\section*{Acknowledgment}
The authors thank Tyrus Berry for insightful discussion and sharing the source codes of his method. The research of JH was partially supported by the Office of Naval Research Grants N00014-13-1-0797, MURI N00014-12-1-0912, and the National Science Foundation grant DMS-1317919. YZ was partially supported as a GRA through NSF grant DMS-1317919.

\section*{Appendix A: Detailed formulation}
In this appendix we review the mathematical formulation for the original Belanger's scheme. In particular, consider the following linear filtering problem:
\BEA
x_{j,k}&=&F_{j,k-1}x_{j,k-1}+\Gamma w_{j,k-1}\label{dsc1}\\
y^o_{j,1}&=&H_{j,1}x_{j,1}+\xi_{j,1}\label{dsc2}
\EEA
where we use the same notations as in Section~2. The only difference is that we use two indices $j,k$, to denote the time point $t_{jN+k}=(jN+k)\delta t$ for $j\geq 0$ and $1\leq k\leq N$, indicating that observations are taken for every $N\geq 1$ integration time steps. Our goal is to estimate the state variable $x_{j,1}$, the covariance $Q$ of the system noise $w$ and the covariance $R$ of the observation noise $\xi$, on the fly. 

Since our goal is to compute the coefficient matrix in \eqref{eq6}, which is a reasonable approximation for \eqref{eq5} when the Kalman gain matrix $K$ is independent of the noises in the previous time steps, let us assume that this condition holds throughout our derivation below. Recall that given observations $y^o_{j,1}$ and prior mean and covariance estimates, $x^f_{j,1}$ and $B^f_{j,1}$ of the state variable $x_{j,1}$ at time $t_{jN+1}$, {\bf Kalman filter} computes the posterior estimate $x^a_{j,1}$ and the posterior covariance matrix $B^a_{j,1}$ of $x_{j,1}$, the prior estimate $x^f_{j,k}$ and the prior covariance $B^f_{j,k}$ of $x_{j,k}$ for $2\leq k\leq N+1$ in the following way:
\BEA
\begin{cases}
K_{j,1}=B^f_{j,1}H_{j,1}^{\top}(\widetilde{R}+H_{j,1}B^f_{j,1}H_{j,1}^{\top})^{-1}\\
x^a_{j,1}=x^f_{j,1}+K_{j,1}(y^o_{j,1}-H_{j,1}x^f_{j,1})\\
B^a_{j,1}=(\mathcal{I}-K_{j,1}H_{j,1})B^f_{j,1}\\
x^f_{j,2}=F_{j,1}x^a_{j,1}\\
B^f_{j,2}=F_{j,1}B^a_{j,1}F_{j,1}^{\top}+\Gamma \widetilde{Q}\Gamma^{\top}\\
x^f_{j,k+1}=F_{j,k}x^f_{j,k} \text{\hspace{15mm}for $k=2,...,N$}\\
B^f_{j,k+1}=F_{j,k}B^a_{j,k}F_{j,k}^{\top}+\Gamma \widetilde{Q}\Gamma^{\top}\text{\hspace{7mm}for $k=2,...,N$}
\end{cases}\label{KF_N}
\EEA
where $\widetilde{Q}$ and $\widetilde{R}$ depends on our knowledge of the true system. If $\widetilde{Q}$ and $\widetilde{R}$ are different from the true value of $Q$ and $R$ the estimates of $x_{j,k}$ would be suboptimal. Define the innovation sequence 
\BEA
v_{j,1}:=y^o_{j,1}-H_{j,1}x^f_{j,1}\label{KF_vk}.
\EEA
Based on equations \eqref{KF_N}, we have:
\BEA
\Delta x^f_{j,1}&:=&x_{j,1}-x^f_{j,1}\nonumber\\
&=&F_{j-1,N}x_{j-1,N}+\Gamma w_{j-1,N}-F_{j-1,N}x^f_{j-1,N}\nonumber\\
&=&F_{j-1,N}F_{j-1,N-1}(x_{j-1,N-1}-x^f_{j-1,N-1})+\Gamma w_{j-1,N}+F_{j-1,N}\Gamma w_{j-1,N-1}\nonumber\\
&=&...\nonumber\\
&=&{\bf F}^{j,1}_{j-1,1}(x_{j-1,1}-x^a_{j-1,1})\nonumber\\
&&+{\bf F}_{j,1}^{j,1}\Gamma w_{j-1,N}+{\bf F}^{j,1}_{j-1,N}\Gamma w_{j-1,N-1}+...+{\bf F}^{j,1}_{j-1,2}\Gamma w_{j-1,1}\label{eq25},
\EEA
where ${\bf F}^{j_1,k_1}_{j_2,k_2}$ denotes the forward operator from time point $t_{j_1N+k_1}$ to time point $t_{j_2N+k_2}$: 
\BEA
{\bf F}_{j_1,k_1}^{j_2,k_2}=\begin{cases}
\mathcal{I} & \text{if $j_1N+k_1=j_2N+k_2$,}\\
\displaystyle\prod_{i=j_1N+k_1}^{j_2N+k_2}F_{i_j,i_k} & \text{if $j_2N+k_2>j_1N+k_1$.}
\end{cases}
\EEA
where $i_j,i_k$ are proper subindexes corresponding to time point $t_i$. We can now write
\eqref{eq25} in a compact form,
\BEA
\Delta x^f_{j,1}&=&\mathcal{U}_{j-1}\Delta x^f_{j-1,1}-\mathcal{G}^{\xi}_{j-1}+\mathcal{F}^w_{j-1},\label{eq27}
\EEA
where 
\BEA
\begin{cases}
\mathcal{U}_{j-1}:={\bf F}^{j,1}_{j-1,1}(I-K_{j-1,1}H_{j-1,1})\\
\mathcal{G}^{\xi}_{j-1}:={\bf F}^{j,1}_{j-1,1}K_{j-1,1}\xi_{j-1,1}\\
\mathcal{F}^w_{j-1}:={\bf F}_{j,1}^{j,1}\Gamma w_{j-1,N}+{\bf F}^{j,1}_{j-1,N}\Gamma w_{j-1,N-1}+...+{\bf F}^{j,1}_{j-1,2}\Gamma w_{j-1,1}
\end{cases}\label{eq28}
\EEA
Let's rewrite \eqref{eq27} in a recursive form,
\BEA
\Delta x^f_{j,1}&=&\mathcal{U}_{j-1}\mathcal{U}_{j-2}\Delta x^f_{j-2,1}-(\mathcal{G}_{j-1}^{\xi}+\mathcal{U}_{j-1}\mathcal{G}_{j-2}^{\xi})+(\mathcal{F}^w_{j-1}+\mathcal{U}_{j-1}\mathcal{F}^w_{j-2})\nonumber\\
&=&...\nonumber\\
&=&(\prod_{i=1}^{j}\mathcal{U}_{j-i})\Delta x^f_{0,1}-\mathcal{G}_{(j)}^{\xi}+\mathcal{F}_{(j)}^w,
\EEA
where 
\BEA
\mathcal{G}^{\xi}_{(j)}&=&\mathcal{G}^{\xi}_{j-1}+\mathcal{U}_{j-1}\mathcal{G}^{\xi}_{j-2}+...+\mathcal{U}_{j-1}...\mathcal{U}_1\mathcal{G}^{\xi}_{0},\label{gxij}\\
\mathcal{F}^{w}_{(j)}&=&\mathcal{F}^{w}_{j-1}+\mathcal{U}_{j-1}\mathcal{F}^{w}_{j-2}+...+\mathcal{U}_{j-1}...\mathcal{U}_1\mathcal{F}^{w}_{0}.\label{fwj}
\EEA
Since $\xi_{j,k}$ and $w_{j,k}$ are i.i.d. standard normal variables and they are also independent of $\Delta x^f_{0,1}$, and since we assume, in this appendix, that Kalman gain matrices $K_{j,k}$ are independent to these noises, $\xi_{j,k}$ and $w_{j,k}$ at time less than $t_{jN+1}$, we would have:
\BEA
\mathbb{E}{\big [}\Delta x^f_{j,1}(\Delta x^f_{j-l,1})^{\top}{\big ]}&=&(\prod_{i=1}^j\mathcal{U}_{j-i})\mathbb{E}{\big [}\Delta x^f_{0,1}(\Delta x^f_{0,1})^{\top}{\big ]}(\prod_{i=1}^{j-l}\mathcal{U}_{j-l-i})^{\top}\nonumber\\
&&+\mathbb{E}{\big [}\mathcal{F}^w_{(j)}(\mathcal{F}^w_{(j-l)})^{\top}{\big ]}+\mathbb{E}{\big [}\mathcal{G}^{\xi}_{(j)}(\mathcal{G}^{\xi}_{(j-l)})^{\top}{\big ]}\label{eq30}
\EEA
where the expectation is taken with respect to each realization of $w$ and $\xi$.
Given that the system is asymptotically stable, $\displaystyle\prod_{i=1}^{j-l}\mathcal{U}_{j-i}$ decays exponentially fast \cite{Belanger1974}. Hence the first term in \eqref{eq30} is small when $j-l$ is large enough. Further it is not hard to see that the last two terms in \eqref{eq30} are linear functions of $Q$ or $R$ respectively, meaning that if we linearly parameterize $Q=\displaystyle\sum_{s=1}^{N_Q}\alpha_sQ_s$ and $R=\displaystyle\sum_{s=1}^{N_R}\beta_sR_s$ using some prescribed basis $Q_s,R_s$, we can write the last two terms of \eqref{eq30} as a linear function of $\alpha_s$ and $\beta_s$:
\BEA
\mathbb{E}{\big [}\mathcal{F}^w_{(j)}(\mathcal{F}^w_{(j-l)})^{\top}{\big ]}&=&\sum_{s=1}^{N_Q}\alpha_s\Phi_{j,l,s}^{(Q)}\label{PhiQ}\\
\mathbb{E}{\big [}\mathcal{G}^{\xi}_{(j)}(\mathcal{G}^{\xi}_{(j-l)})^{\top}{\big ]}&=&\sum_{s=1}^{N_R}\beta_s\Phi_{j,l,s}^{(R)}.\label{PhiR}
\EEA
Belanger pointed out that the matrix $\Phi^{(Q)}_{j,l,s}$ and $\Phi^{(R)}_{j,l,s}$ can be computed recursively \cite{Belanger1974}:
\BEA\label{PsiQ}
\Phi^{(Q)}_{j,l,s}&=&\begin{cases}
\mathcal{U}_{j-1}\Phi^{(Q)}_{j-1,l-1,s} & \text{if $l>0$}\\
&\\
\mathcal{U}_{j-1}\Phi^{(Q)}_{j-1,0,s}\mathcal{U}_{j-1}^{\top}+\Gamma Q_s\Gamma^{\top} & \text{if $l=0$, $N=1$}\\
&\\
\!\begin{aligned}
& \mathcal{U}_{j-1}\Phi^{(Q)}_{j-1,0,s}\mathcal{U}_{j-1}^{\top}+\Gamma Q_s\Gamma^{\top}\\
& +\displaystyle\sum_{i=1}^{N-1}{\bf F}^{j,1}_{j-1,i+1}\Gamma Q_s\Gamma^{\top}({\bf F}_{j-1,i+1}^{j,1})^{\top} \end{aligned} & \text{if $l=0$, $N>1$}
\end{cases}\\ \label{PsiR} 
\Phi^{(R)}_{j,l,s}&=&\begin{cases}
\mathcal{U}_{j-1}\Phi^{(R)}_{j-1,l-1,s} & \text{if $l>0$}\\
\mathcal{U}_{j-1}\Phi^{(R)}_{j-1,0,s}\mathcal{U}_{j-1}^{\top}+\mathcal{S}_{j-1}R_s\mathcal{S}_{j-1}^{\top} & \text{if $l=0$},
\end{cases}
\EEA
where we define $\mathcal{S}_{j-1}:={\bf F}^{j,1}_{j-1,1}K_{j-1,1}$. Based on equations \eqref{KF_vk}, \eqref{eq30}-\eqref{PsiR}, we have
\BEA
\mathbb{E}(v_{j,1}v_{j-l,1}^{\top})&=&H_{j,1}\mathbb{E}{\big [}\Delta x^f_{j,1}(\Delta x^f_{j-l,1})^{\top}{\big ]}H_{j-l,1}^{\top}+H_{j,1}\mathbb{E}{\big [}\Delta x^f_{j,1}\xi_{j-l,1}^{\top}{\big ]}\nonumber\\
&=&\sum_{s=1}^{N_Q}\alpha_sH_{j,1}\Phi^{(Q)}_{j,l,s}H_{j-l,1}^{\top}+\sum_{s=1}^{N_R}\beta_sH_{j,1}\Phi^{(R)}_{j,l,s}H_{j-l,1}^{\top}\nonumber\\
&&+\begin{cases}
\sum_{s=1}^{N_R}\beta_sR_s & \text{if $l=0$}\\
-\sum_{s=1}^{N_R}\beta_sH_{j,1}\mathcal{S}_{j-1}R_s & \text{if $l=1$}\\
-\sum_{s=1}^{N_R}\beta_sH_{j,1}\mathcal{U}_{j-1}...\mathcal{U}_{j-l+1}\mathcal{S}_{j-l}R_s & \text{if $l>1$.}
\end{cases}\label{eq35}
\EEA
and we can define a linear observation operator for $Q$ and $R$ with coefficients:
\BEA
\mathcal{H}^{(Q)}_{j,l,s}&:=&H_{j,1}\Phi^{(Q)}_{j,l,s}H_{j-l,1}^{\top} \label{HQjls}\\
\mathcal{H}^{(R)}_{j,l,s}&:=&H_{j,1}\Phi^{(R)}_{j,l,s}H_{j-l,1}^{\top}+\begin{cases}
R_s & \text{if $l=0$}\\
-H_{j,1}\mathcal{S}_{j-1}R_s & \text{if $l=1$}\\
-H_{j,1}\mathcal{U}_{j-1}...\mathcal{U}_{j-l+1}\mathcal{S}_{j-l}R_s & \text{if $l>1$},\label{HRjls}
\end{cases}
\EEA
so that \eqref{eq35} can be conveniently written as
\BEA
\mathbb{E}(v_{j,1}v_{j-l,1}^{\top})=\sum_{s=1}^{N_Q}\alpha_s\mathcal{H}^{(Q)}_{j,l,s}+\sum_{s=1}^{N_R}\beta_s\mathcal{H}^{(R)}_{j,l,s},\label{eq38}
\EEA
which is exactly \eqref{eq6} in the case when observations are taken every integration time step ($N=1$).

\section*{Appendix B. Pseudo-code}
In this appendix, we present a complete pseudo-code of the modified Belanger's method. To make this paper self-contained we also present the ETKF algorithm that we used for filtering stochastic nonlinear systems. This version of ETKF is similar to the one used by  \cite{HarlimMahdiMajda2014} and it is a variant of the one formulated earlier by \cite{hunt:07}. We should also mention that this pseudo-code is not optimal so one may need to optimize it carefully for high-dimensional problems. The following notations are needed:
\begin{nolinenumbers}
{\scriptsize
\begin{itemize}
\setlength{\columnsep}{20pt}
\begin{multicols}{2}
\setlength\itemsep{0.1em}
\item $n\in\mathbb{N}$: the dimension of state variable;
\item $m\in\mathbb{N}$: the dimension of observation;
\item $N\in \mathbb{N}$: number of time steps per observation;
\item $N_e\in\mathbb{N}$: the ensemble size;
\item $n_s\in\mathbb{N}$: the dimension of system noise $w_k$;
\item $x_{j,k}\in\mathbb{R}^{n\times 1}$: the true state variable at the  $jN+k$-th  time step;
\item $X^{f,i}_{j,k}\in\mathbb{R}^{n\times 1}$:   the $i$-th ensemble member of prior estimate of state variable at the  $jN+k$-th  time step, for $i=1,...,N_e$;
\item $X^{df,i}_{j,k}\in\mathbb{R}^{n\times 1}$:   the $i$-th ensemble member of deterministic forecast  of state variable at the  $jN+k$-th time step, for $i=1,...,N_e$;
\item $\bar{X}^{f}_{j,k}\in\mathbb{R}^{n\times 1}$:  ensemble mean of prior estimate of state variable at the  $jN+k$-th  time step;
\item $X^{a,i}_{j,k}\in\mathbb{R}^{n\times 1}$: the $i$-th ensemble member of posterior estimate of state variable at the  $jN+k$-th  time step, for $i=1,...,N_e$;
\item $\bar{X}^a_{j,k}\in\mathbb{R}^{n\times 1}$: ensemble mean of posterior estimate of state variable at the  $jN+k$-th  time step;
\item $v_{j,1}\in\mathbb{R}^{m\times 1}$: the innovation at the  $jN+1'$st  time step;
\item $B^f_{j,k}\in\mathbb{R}^{n\times n}$: prior covariance of state variable at the  $jN+k$-th  time step;
\item $B^a_{j,1}\in\mathbb{R}^{n\times n}$: posterior covariance of the state variable at the  $jN+1'$st  time step;
\item $K_{j,k}\in\mathbb{R}^{n\times m}$: the Kalman gain at the  $jN+k$-th  time step;
\item $h_{j,1}$:  observational operator at the $jN+1'$st  time step;
\item $H_{j,1}\in\mathbb{R}^{m\times n}$:  the linear version of $h_{j,1}$;
\item $f_{j,k}$: the nonlinear deterministic model operator that transmits state variable from the $jN+k$-th time step to the $jN+k+1'$st time step;
\item $F_{j,k}\in\mathbb{R}^{n\times n}$: the linear version of $f_{j,k}$;
\item $Q\in\mathbb{R}^{n_s\times n_s}$: the covariance of system noise $w_{k}$;
\item $\Gamma\in\mathbb{R}^{n\times n_s}$: the coefficient of system noise;
\item $R\in\mathbb{R}^{m\times m}$: the covariance of observation noise $\xi_k$.
\item $\Phi^{(Q)}_{j,l,s}$ and  $\Phi^{(R)}_{j,l,s}$ : the matrices defined by \eqref{PhiQ} and \eqref{PhiR} for $l=0,...,L$;
\item $\mathcal{H}^{(Q)}_{j,l,s}$ and  $\mathcal{H}^{(R)}_{j,l,s}$ : the matrices defined by \eqref{eq38} for $l=0,...,L$;
\item $\mathcal{H}^{(Q)\text{sum}}_{J,l,s}$ and  $\mathcal{H}^{(R)\text{sum}}_{J,l,s}$ : the sum of $\mathcal{H}^{(Q)}_{j,l,s}$ and  $\mathcal{H}^{(R)}_{j,l,s}$, respectively,  from time $t_{(L+1)N+1}$ to $t_{JN+1}$ 
\item $\mathcal{Y}^{(QR){\text{sum}}}_{J,l}$ : the sum of $v_{j,1}v_{j-l,1}^{\top}$ from $t_{(L+1)N+1}$ up to time $t_{JN+1}$;
\item $N_Q$ : number of parameters for $Q$;
\item $N_R$ : number of parameters for $R$;
\item $\mathcal{Q}_s$ : the basis of parameterization of $Q$, for $s=1,...,N_Q$;
\item $\mathcal{R}_s$ : the basis of parameterization of $R$, for $s=1,...,N_R$;
\item $Q_j$ : the estimation of $Q$ that is used in the primary filter at the $jN+1'$st time step;
\item $R_j$ : the estimation of $R$ that is used in the primary filter at the $jN+1'$st time step;
\item $\tau$ : the relaxation coefficient in \eqref{Relx1} and \eqref{Relx2}.
\end{multicols}
\end{itemize}}
\end{nolinenumbers}

\begin{algorithm}[H]
\caption{Covariance Estimation at the $jN+1$'s Time Step} 
\label{CovEst_k} 
\begin{algorithmic}[1] 
\REQUIRE 
$v_{j-l,1},B^f_{j-l,1},  K_{j-l,1}, F_{j-l-1,k}, H_{j-l,1}, \mathcal{Q}_s,\mathcal{R}_s, \Gamma, \alpha_{s,j-1}, Q_{j-1}, R_{j-1}, \beta_{s,j-1}$, $\mathcal{U}_{j-l-1},  \mathcal{S}_{j-l-1}, \Phi^{(Q)}_{j-1,l,s}, \Phi^{(R)}_{j-1,l,s},  \mathcal{H}^{(Q){\text{sum}}}_{j-1,l,s}, \mathcal{H}^{(R){\text{sum}}}_{j-1,l,s}, \mathcal{Y}^{(QR){\text{sum}}}_{j-1,l}, \tau$, $\text{for } 1\leq s\leq N_Q \text{ or } N_R, 0\leq l\leq L$
\ENSURE 
$\alpha_{s,j}, \beta_{s,j}, Q_j, R_j, \Phi^{(Q)}_{j,l,s},  \Phi^{(R)}_{j,l,s}, \mathcal{Y}^{(QR){\text{sum}}}_{j,l}, \mathcal{H}^{(Q){\text{sum}}}_{j,l,s}, \mathcal{H}^{(R){\text{sum}}}_{j,l,s}$, for $0\leq l\leq L$
 
 \IF{$j<L+1$} \mycommenta{do not do anything for the first $L$ time steps }
 \STATE $\alpha_{s,j}\gets \alpha_{s,j-1}$\hspace{5mm} for $s=1,...,N_Q$
 \STATE $\beta_{s,j}\gets \beta_{s,j-1}$\hspace{5mm} for $s=1,...,N_R$
 \STATE $Q_j\gets Q_{j-1}$
 \STATE $R_j\gets R_{j-1}$
 \STATE $\Phi^{(Q)}_{j,l,s}\gets \Phi^{(Q)}_{j-1,l,s}$\hspace{5mm} for $l=0,...,L$ and $s=1,...,N_Q$
 \STATE $\Phi^{(R)}_{j,l,s}\gets \Phi^{(R)}_{j-1,l,s}$\hspace{5mm} for $l=0,...,L$ and $s=1,...,N_R$
  \STATE $\mathcal{Y}^{(QR){\text{sum}}}_{j,l}\gets \mathcal{Y}^{(QR){\text{sum}}}_{j-1,l}$\hspace{5mm} for $l=0,...,L$  
  \STATE $\mathcal{H}^{(Q){\text{sum}}}_{j,l,s}\gets \mathcal{H}^{(Q){\text{sum}}}_{j-1,l,s}$\hspace{5mm} for $l=0,...,L$ and $s=1,...,N_Q$
 \STATE $\mathcal{H}^{(R){\text{sum}}}_{j,l,s}\gets \mathcal{H}^{(R){\text{sum}}}_{j-1,l,s}$\hspace{5mm} for $l=0,...,L$ and $s=1,...,N_R$
 \ELSE
 \FOR{$l=0$ to $L$}
\STATE \fbox{$\mathcal{Y}^{(QR){\text{sum}}}_{j,l}\gets \mathcal{Y}^{(QR){\text{sum}}}_{j-1,l}+v_{j,1}v_{j-l,1}^{\top}$} \mycommenta{Compute $\displaystyle\sum_{j=L+1}^{J}v_jv_{j-l}^{\top}$ in \eqref{newmethod}}
\IF{$l=0$}
\FOR{$s=1$ to $N_Q$}
\STATE \fbox{$\Phi^{(Q)}_{j,l,s}\gets \mathcal{U}_{j-1}\Phi^{(Q)}_{j-1,0,s}\mathcal{U}_{j-1}^{\top}+\mathcal{Q}_s$}
\STATE ${\bf F}^{j,1}_{j,1}\gets I$
\algstore{CovEst1}
 \end{algorithmic}
\end{algorithm}

\begin{algorithm}                     
\begin{algorithmic} [H]                   
\algrestore{CovEst1}
\FOR{$k=N$ to $2$}
\STATE ${\bf F}^{j,1}_{j-1,k}\gets {\bf F}^{j,1}_{j-1,k+1}F_{j-1,k}$
\STATE \fbox{$\Phi^{(Q)}_{j,1,s}\gets \Phi^{(Q)}_{j-1,1,s}+{\bf F}^{j,1}_{j-1,k}\Gamma\mathcal{Q}_s\Gamma^{\top}({\bf F}^{j,1}_{j-1,k})^{\top}$}
\ENDFOR
\STATE $\mathcal{H}^{(Q)}_{{\text{temp}},s}\gets H_{j,1}\Phi^{(Q)}_{j,1,s}H_{j,1}^{\top}$\mycommenta{Compute $\mathcal{H}^{(Q)}_{j,l,s}$ in \eqref{newmethod}}
\ENDFOR
\FOR{$s=1$ to $N_R$}
\STATE \fbox{$\Phi^{(R)}_{j,1,s}\gets \mathcal{U}_{j-1}\Phi^{(R)}_{j,1,s}\mathcal{U}_{j-1}^{\top}+\mathcal{S}_{j-1}\mathcal{R}_{s}\mathcal{S}_{j-1}^{\top}$}
\STATE $\mathcal{H}^{(R)}_{\text{temp},s}\gets H_{j,1}\Phi^{(R)}_{j,1,s}H_{j,1,s}^{\top}+\mathcal{R}_s$
\ENDFOR
\ELSE
\FOR{$s=1$ to $N_Q$}
\STATE \fbox{$\Phi^{(Q)}_{j,l,s}\gets \mathcal{U}_{j-1}\Phi^{(Q)}_{j-1,l-1,s}$}
\STATE $\mathcal{H}^{(Q)}_{\text{temp},s}\gets H_{j,1}\Phi^{(Q)}_{j,l,s}H_{j,1}^{\top}$
\ENDFOR
\FOR{$s=1$ to $N_R$}
\STATE \fbox{$\Phi^{(R)}_{j,1,s}\gets \mathcal{U}_{j-1}\Phi^{(R)}_{j-1,l-1,s}$}
\IF{$l=1$}\mycommenta{Compute $\mathcal{H}^{(R)}_{j,l,s}$ in \eqref{newmethod}}
\STATE $\mathcal{H}^{(R)}_{\text{temp},s}\gets H_{j,1}\Phi^{(R)}_{j,l,s}H_{j-l,1}^{\top}-H_{j,1}\mathcal{S}_{j-1}\mathcal{R}_s$
\ELSE
\STATE $\mathcal{H}^{(R)}_{\text{temp},s}\gets H_{j,1}\Phi^{(R)}_{j,l,s}H_{j-l,s}^{\top}-H_{j,1}\mathcal{U}_{j-1}...\mathcal{U}_{j-l+1}\mathcal{S}_{j-1}\mathcal{R}_s$
\ENDIF
\ENDFOR
\ENDIF
\ENDFOR
\algstore{CovEst2}
 \end{algorithmic}
\end{algorithm}

\begin{algorithm}                     
\begin{algorithmic} [H]                   
\algrestore{CovEst2}
\STATE \mycommentb{Compute $\displaystyle\sum_{j=L+1}^J\mathcal{H}^{(Q)}_{j,l,s}$ and $\displaystyle\sum_{j=L+1}^J\mathcal{H}^{(R)}_{j,l,s}$ in \eqref{newmethod} }
\STATE \fbox{$\mathcal{H}^{(Q){\text{sum}}}_{j,l,s}\gets \mathcal{H}^{(Q){\text{sum}}}_{j-1,l,s}+\mathcal{H}^{(Q)}_{\text{temp},s}$}, \hspace{5mm} for $l=0,1,...,L$ and $s=1,...,N_Q$
\STATE \fbox{$\mathcal{H}^{(R){\text{sum}}}_{j,l,s}\gets \mathcal{H}^{(R){\text{sum}}}_{j-1,l,s}+\mathcal{H}^{(R)}_{\text{temp},s}$}, \hspace{5mm} for $l=0,1,...,L$ and $s=1,...,N_R$
\STATE reshape $\mathcal{H}^{(Q){\text{sum}}}_{j,l,s}$ to get $\mathfrak{h}^{(q)}_{l,s}\in \mathbb{R}^{m^2\times 1}$ for $1\leq s \leq N_Q$ and $0\leq l\leq L$
\STATE reshape $\mathcal{H}^{(R){\text{sum}}}_{j,l,s}$ to get $\mathfrak{h}^{(r)}_{l,s}\in \mathbb{R}^{m^2\times 1}$ for $1\leq s\leq N_R$ and $0\leq l\leq L$
\STATE $\mathfrak{H}^{(q)}_l\gets [\mathfrak{h}^{(q)}_{l,1},\mathfrak{h}^{(q)}_{l,2},...,\mathfrak{h}^{(q)}_{l,N_Q}]$ \hspace{5mm} for $l=0,1,...,L$
\STATE $\mathfrak{H}^{(r)}_l\gets [\mathfrak{h}^{(r)}_{l,1},\mathfrak{h}^{(r)}_{l,2},...,\mathfrak{h}^{(r)}_{l,N_R}]$ \hspace{5mm} for $l=0,1,...,L$
\STATE $\mathfrak{H}^{(q)}\gets \begin{pmatrix}
\mathfrak{h}^{(q)}_{0} \\
\mathfrak{h}^{(q)}_{1}\\
\vdots\\
\mathfrak{h}^{(q)}_{L}
\end{pmatrix}$, \hspace{5mm}$\mathfrak{H}^r\gets \begin{pmatrix}
\mathfrak{h}^{(r)}_{0} \\
\mathfrak{h}^{(r)}_{1}\\
\vdots\\
\mathfrak{h}^{(r)}_{L}
\end{pmatrix}, \hspace{5mm}
\mathfrak{Y}^{{(qr)}}\gets \begin{pmatrix}
\mathcal{Y}^{(QR){\text{sum}}}_{j,0} \\
\mathcal{Y}^{(QR){\text{sum}}}_{j,1}\\
\vdots\\
\mathcal{Y}^{(QR){\text{sum}}}_{j,L}
\end{pmatrix}$
\STATE $\mathfrak{H}^{(qr)}\gets [\mathfrak{H}^{(q)},\mathfrak{H}^{(r)}]$
\STATE $\Lambda^{(qr)}\gets ((\mathfrak{H}^{(qr)})^{\top}\mathfrak{H}^{(qr)})^{\dagger}(\mathfrak{H}^{(qr)})^{\top}\mathfrak{Y}^{(qr)}$, where $\dagger$ denotes the matrix pseudo-inverse. \mycommenta{Solve \eqref{newmethod} for the covariance parameters using least square method}
\STATE \mycommentb{Temporal smoothing of covariance parameters. Refer to \eqref{Relx1} and \eqref{Relx2}}
\STATE \fbox{$\alpha_{s,j}\gets \alpha_{s,j-1}+1/\tau(\Lambda^{(qr)}_s-\alpha_{s,j-1})$} \hspace{5mm} for $s=1,...,N_Q$
\STATE  \fbox{$\beta_{s,j}\gets \beta_{s,j-1}+1/\tau(\Lambda^{(qr)}_{s+N_Q}-\beta_{s,j-1})$} \hspace{5mm} for $s=1,...,N_R$
\STATE \fbox{$Q_{j}\gets \displaystyle\sum_{s=1}^{N_Q}\alpha_{s,j}\mathcal{Q}_{s}$},\fbox{ $R_{j}\gets \displaystyle\sum_{s=1}^{N_R}\beta_{s,j}\mathcal{R}_{s}$ }\mycommenta{Construct $Q$ and $R$ from the newly estimated parameters}
 \ENDIF     
\end{algorithmic}
\end{algorithm}

\begin{algorithm} 
\caption{Ensemble Transform Kalman Filter (for stochastic nonlinear systems) from the $(jN+1)'$st Time Step to $(j+1)N+1'$st Time Step} 
\label{ETKF_k} 
\begin{algorithmic}[1] 
\REQUIRE 
$y^o_{j,1}, X^{f,i}_{j,1}, B^f_{j,1},  f_{j,k}, h_{j,1}, Q, R, \Gamma$

\ENSURE 
$X^{f,i}_{j+1,1}, B^f_{j+1,1}, F_{j,k},H_{j,1}, K_{j,1}, B^a_{j,1}$ for $1\leq k\leq N$

\STATE \mycommentb{ETKF is widely known and has many different versions. For a detailed description of this version please refer to \cite{hunt:07, HarlimMahdiMajda2013} }
\STATE $\bar{X}^f_{j,1}\gets \frac{1}{N_e}\displaystyle\sum_{i=1}^{N_e}X^{f,i}_{j,1}$
\STATE $U^{(i)}_{j,1}\gets X^{f,i}_{j,1}-\bar{X}^{f}_{j,1}$, \hspace{5mm} for $i=1,2,...,N_e$.
\STATE $Y^{f,i}_{j,1}\gets h_{j,1}(X^{f,i}_{j,1})$, \hspace{5mm} for $i=1,2,...,N_e$.
\STATE $U\gets [U^{f,1}_{j,1},...,U^{f,N_e}_{j,N_e}]$
\STATE $\bar{Y}^{f}_{j,1}\gets \frac{1}{N_e}\displaystyle\sum_{i=1}^{N_e}Y^{f,i}_{j,1}$
\STATE $V^{f,i}_{j,1}\gets Y^{f,i}_{j,1}-\bar{Y}^{f}_{j,1}$, \hspace{5mm} for $i=1,2,...,N_e$
\STATE $V\gets [V^{f,1}_{j,1},...,V^{f,N_e}_{j,1}]$
\STATE \fbox{$H_{j,1}\gets VU^{\dagger}$} where $\dagger$ refers to the matrix pseudo-inverse.
\STATE $J\gets (N_e-1)I+V^{\top}R^{-1}V$
\STATE $x\gets J^{-1}V^{\top}R^{-1}v_{j,1}$
\STATE $\bar{X}^a_{j,1}\gets \bar{X}^{f}_{j,1}+Ux$
\STATE Compute the singular value decomposition (SVD) of $J=E_1D_1E_1^{\top}$
\STATE $D_2\gets (N_e-1)D_1^{-1}$     
\STATE $T\gets E_1\sqrt{D_2}E_1^{\top}$
\STATE $P_y\gets \frac{1}{N_e-1}VV^{\top}+R$
\STATE $P_{xy}\gets \frac{1}{N_e-1}UV^{\top}$
\STATE \fbox{$K_{j,1}\gets P_{xy}P_y^{-1}$}
\STATE $U^{a}\gets UT $, denote the $i$-th column of $U^a$ by $U^{a,i}$;
\STATE \fbox{$B^a_{j,1}\gets U^a(U^a)^{\top}/(En-1)$}
\STATE $X^{a,i}_{j,1}\gets \bar{X}^a_{j,1}+U^{a,i}$,\hspace{5mm} for $i=1,2,...,N_e$.
    \algstore{ETKF_k}
    \end{algorithmic}
\end{algorithm}

\begin{algorithm}                     
\begin{algorithmic} [1]                   
\algrestore{ETKF_k}
\STATE $X^{\text{temp}}\gets X^a_{j,1}$, denote the $i$-th column by $X^{\text{temp},i}$
\FOR{$k=1$ to $N$}
\STATE $\bar{X}^{\text{temp}}\gets \frac{1}{N_e}\displaystyle\sum_{i=1}^{N_e}X^{\text{temp},i}$
\STATE $U^{f,i}\gets X^{\text{temp},i}-\bar{X}^{\text{temp}}$\hspace{5mm} for $i=1,2,...,N_e$.
\STATE $U^f\gets [U^{f,1},...,U^{f,N_e}]$
\STATE $B^{\text{temp}}\gets \frac{1}{N_e-1}U^f(U^f)^{\top}$
\STATE $X^{\text{temp},i}\gets f_{j,k}(X^{\text{temp},i})$
\STATE $\bar{X}^{\text{temp}}\gets \frac{1}{N_e}\displaystyle\sum_{i=1}^{N_e}X^{\text{temp},i}$
\STATE $U^{df,i}\gets X^{\text{temp},i}-\bar{X}^{\text{temp}}$\hspace{5mm} for $i=1,2,...,N_e$
\STATE $U^{df}\gets [U^{df,1},...,U^{df,N_e}]$
\STATE \fbox{$F_{j,k}\gets U^{df}(U^f)^{\dagger}$}
\STATE $B^{\text{temp}}\gets \frac{1}{N_e-1}U^{f}(U^{f})^{\top}+\Gamma Q\Gamma^{\top}$
\STATE Independently generate $N_e$ standard normal random vectors $\delta^i\in\mathbb{R}^{n\times 1}$, for $1\leq i\leq N_e$
\STATE $\bar{\delta}\gets \frac{1}{N_e}\displaystyle\sum_{i=1}^{N_e}\delta^i$
\STATE $\delta^i\gets \delta^i-\bar{\delta}$\hspace{5mm} for $i=1,2,...,N_e$.
\STATE $\delta\gets[\delta^1,...,\delta^{N_e}]$
\STATE Compute the singular value decomposition of $\frac{1}{N_e-1}\delta\delta^{\top}=E_3D_3E_3^{\top}$
\STATE $U^{\text{temp}}\gets \sqrt{B^{\text{temp}}}E_3(D_3)^{-1}E_3^{\top}\delta$, denote the $i$-th column by $U^{\text{temp},i}$
\STATE $X^{\text{temp},i}\gets \bar{X}^{\text{temp}}+U^{\text{temp},i}$\hspace{5mm} for $i=1,2,...,N_e$
\ENDFOR
\STATE \fbox{$X^{f,i}_{j+1,1}\gets X^{\text{temp},i}$ for $i=1,...,N_e$}
\STATE \fbox{$B^f_{j+1,1}\gets B^{\text{temp}}$}
\end{algorithmic}
\end{algorithm}

\clearpage

\end{document}